\documentclass[brochure,12pt]{bourbaki} 
\usepackage[matrix,arrow]{xy}
\usepackage{amssymb,amsfonts,amsmath,footnote}
\usepackage[francais]{babel}
\date{14 Novembre 2009}
\bbkannee{62\`eme ann\'ee, 2009-2010}
\bbknumero{1012}
\title{
GROUPE DE CHOW DES Z\'ERO-CYCLES\\ SUR LES VARI\'ET\'ES $p$-ADIQUES }
\subtitle{d'apr\`es  S. Saito, K. Sato
 et  al.}
\author{Jean-Louis COLLIOT-TH\'EL\`ENE}
\address{
CNRS\\Universit\'e Paris XI\\
UMR 8628 du CNRS\\
D\'epartement de Math\'ematiques\\
B\^atiment 425\\
F--91405 ORSAY C\'edex}
\email{jlct@math.u-psud.fr}

\newcommand{\et}{{\textup{\'et}}}

\def\cal{\mathcal}

\def\N{{\bf N}}
\def\F{{\bf F}}
\def\C{{\bf C}}
\def\Q{{\bf Q}}
\def\Z{{\bf Z}}
\def\G{{\bf G}}
\def\K{{\mathcal K}}
\def\br{{\rm Br} } 
\def\pic{{\rm Pic}}
\def\P{{\bf P}}
\def\spec{{\rm Spec}\ }
\def\X{{\mathcal X}}
\def\deg{{\rm deg}}
\def\div{{\rm div}}
\def\NS{{\rm NS}}
\def\alb{{\rm alb}}
\def\Alb{{\rm Alb}}
\def\Hom{{\rm Hom}}
\def\dim{{\rm dim}}

\def\cqfd{\hfill $\square$}
 
\def\oi{\hskip1mm {\buildrel \simeq \over \rightarrow} \hskip1mm}
\def\loi{\hskip1mm {\buildrel \simeq \over \longrightarrow} \hskip1mm}

\def\lloi{\hskip1mm {\buildrel \simeq \over \longleftarrow} \hskip1mm}

\newcommand{\colim}{\operatornamewithlimits\varinjlim}

\newcommand{\colime}{\mathop{\lim}\limits_{{\longrightarrow}\atop n}}

\usepackage[T1]{fontenc} 
\begin{document}
\maketitle

 \bigskip

\bigskip
 
  \noindent{\bf    INTRODUCTION}

\bigskip

Soient $k$ un corps et $X$ une $k$-vari\'et\'e alg\'ebrique projective, lisse et g\'eom\'etriquement
irr\'eductible (cette derni\`ere hypoth\`ese sera souvent tacitement faite). 
 On note $Z_{0}(X)$  le groupe des z\'ero-cycles sur $X$, c'est-\`a-dire le groupe ab\'elien libre sur les points ferm\'es de $X$ (points $x$ du sch\'ema $X$ dont le corps r\'esiduel $\kappa(x)$ est une extension finie de $k$). On dispose d'une application degr\'e
$$ \deg : Z_{0}(X) \to \Z$$
d\'efinie par lin\'earit\'e \`a partir de l'application envoyant un point ferm\'e $x$ sur 
le degr\'e $[\kappa(x):k]$.

\`A tout couple form\'e d'une courbe ferm\'ee int\`egre $C \subset X$ et d'une fonction rationnelle
non nulle $f \in k(C)^{\times}$, on associe un z\'ero-cycle, le diviseur de $f$. Celui-ci est ainsi d\'efini :
on consid\`ere la normalisation $\tilde{C} \to C$ de la courbe $C$, et le morphisme compos\'e
$\pi : \tilde{C} \to C \to X$. On d\'efinit alors $\div(f) = \pi_{*}(\div_{\tilde{C}}(f))$.
Le groupe de Chow des z\'ero-cycles sur $X$ est par d\'efinition le quotient de $Z_{0}(X)$
par le sous-groupe engendr\'e par tous les $\div(f)$ pour tous les couples $(C,f)$.

Comme la $k$-vari\'et\'e $X$ est projective,  l'application degr\'e induit un homomorphisme
$\deg : CH_{0}(X) \to \Z$. On note  $A_{0}(X)$ le noyau de cette application. On dispose
donc d'une suite exacte
$$0 \to A_{0}(X) \to CH_{0}(X) \to \Z.$$
L'image de la fl\`eche degr\'e est un sous-groupe $\Z.I_{X} \subset \Z$ d'indice fini.
La suite est scind\'ee si $X$ poss\`ede un z\'ero-cycle 
$z_{0}$ de degr\'e 1, {\it ce qu'on suppose 
d\'esormais dans cette introduction}.

\`A la $k$-vari\'et\'e $X$ on associe  sa vari\'et\'e de Picard $\pic_{X/k, red}^0$, qui
est
une vari\'et\'e
ab\'elienne.  
On a une suite exacte
$$0 \to \pic_{X/k,red}^0(k) \to \pic X \to \NS(X) \to 0$$
o\`u $\NS(X)$ est un groupe ab\'elien de type fini.

La vari\'et\'e ab\'elienne duale de $\pic_{X/k,red}^0$ est la vari\'et\'e d'Albanese $ \Alb_{X/k}$ de $X$.
\`A la donn\'ee de $z_{0}$ est associ\'e un $k$-morphisme  $$\alb_{X} :X \to  \Alb_{X/k}$$
 induisant un isomorphisme sur les vari\'et\'es de Picard de ces deux vari\'et\'es. 
 Ce morphisme induit un homomorphisme de groupes ab\'eliens $$\alb_{X} : A_{0}(X) \to \Alb_{X/k}(k)$$
 qui ne d\'epend pas du choix de $z_{0}$.

\medskip

Lorsque $\dim(X)=1$, c'est-\`a-dire lorsque $X$ est une courbe (projective, lisse)
on a un isomorphisme $\pic X \oi CH_{0}(X)$, qui induit un isomorphisme
$\pic_{X/k,red}^0(k) \oi A_{0}(X)$. La fl\`eche $\alb_{X} : A_{0}(X) \to \Alb_{X/k}(k)$
est un isomorphisme.

Les propri\'et\'es
des groupes de points rationnels de vari\'et\'es ab\'eliennes
donnent alors des th\'eor\`emes sur la structure des groupes $CH_{0}(X)$
et $A_{0}(X)$. En particulier, pour $X/k$ une courbe de genre $g$ avec ou sans z\'ero-cycle de degr\'e 1,
on a les propri\'et\'es suivantes:

\medskip

(1) Si  $k$ est  un corps de type fini sur le corps premier, 
 le groupe $CH_{0}(X)$ est un groupe ab\'elien de type fini (Mordell-Weil).
 
 \medskip

(2) Si $k$ est  un corps fini,  le groupe   $A_{0}(X)$ est fini.

\medskip
(3) Si $k$ est un corps $p$-adique  (ce qui dans cet expos\'e 
signifie extension finie du corps $p$-adique $\Q_{p}$),
 le groupe  $A_{0}(X)$ est extension d'un groupe fini
par un sous-groupe isomorphe \`a une somme directe de $g$
exemplaires de l'anneau des entiers de $k$ (Lutz, Mattuck).

  En cons\'equence, 

(3.1) Le groupe $A_{0}(X)$ est somme directe
d'un groupe fini (d'ordre premier \`a $p$)  et
d'un groupe $p'$-divisible (c'est-\`a-dire 
divisible par tout entier premier \`a $p$).

(3.2) Pour  tout entier  $n>0$, le quotient $CH_{0}(X)/n$ est fini.

(3.3) Pour presque tout premier $l$, on a $A_{0}(X)/l=0$.

(3.4) Le sous-groupe de torsion de $CH_{0}(X)$ est fini.

On peut  en outre d\'etecter
les classes dans $CH_{0}(X)$  au moyen de la cohomologie
\'etale sur $X$ (voir le paragraphe 1 ci-apr\`es).

\medskip

Il est naturel de se demander si certaines parmi ces  propri\'et\'es du groupe $CH_{0}(X)$ valent encore
pour une $k$-vari\'et\'e projective lisse $X$ de dimension quelconque.

Dans la situation  (1),  m\^eme pour $k$ le  corps des rationnels,
en dehors des cas qui se r\'eduisent formellement au th\'eor\`eme de Mordell-Weil,
on n'a aucun r\'esultat non trivial sur la finitude de $CH_{0}(X)/n$ pour $n>1$
ou sur la finitude de la dimension du \mbox{$\Q$-vectoriel} $CH_{0}(X) \otimes_{\Z}\Q$.

Dans la situation  (2), qui porte sur le cas des corps finis,
la finitude de $A_{0}(X)$ est un th\'eor\`eme de K. Kato et S. Saito \cite{KS}.
On en sait beaucoup plus : voir \`a ce sujet l'expos\'e r\'ecent de T. Szamuely  \cite{Sz} sur le corps de classes
de dimension sup\'erieure.

 Le pr\'esent expos\'e porte sur le cas des corps $p$-adiques.

  Depuis les ann\'ees 1980, une m\'ethode de $K$-th\'eorie alg\'ebrique
invent\'ee par S. Bloch  et reposant sur  un  th\'eor\`eme de Merkur'ev et Suslin
  a permis d'obtenir  un certain nombre de r\'esultats, en particulier pour les surfaces.
On \'evoquera ces r\'esultats
  au paragraphe~2.
  
En 2006, S. Saito et K. Sato \cite{SS2}  r\'ealis\`erent que pour obtenir des
\'enonc\'es g\'en\'eraux il vaut mieux consid\'erer non le groupe de Chow  des  z\'ero-cycles sur
une vari\'et\'e projective  et lisse $X$ sur un corps $p$-adique $k$, mais le groupe de
Chow des 1-cycles sur un mod\`ele r\'egulier et projectif de $X$ au-dessus de l'anneau des entiers
de $k$ (lorsqu'un tel mod\`ele existe).  

On d\'ecrira en d\'etail leur travail au paragraphe 3. Le th\'eor\`eme principal est
le th\'eor\`eme \ref{thmprincipal}. En voici deux applications
 (Th\'eor\`emes \ref{bonnereduc} et  \ref{applicgeneral}).

\begin{theo} 
Soit $X$ une vari\'et\'e projective, lisse, g\'eom\'etriquement connexe sur un corps $p$-adique $k$.
Si $X$ a bonne r\'eduction $Y$ sur le corps r\'esiduel fini $F$, alors la fl\`eche de sp\'ecialisation
$A_{0}(X) \to A_{0}(Y)$, qui est une surjection sur le groupe fini $ A_{0}(Y)$, a
un noyau $p'$-divisible, c'est-\`a-dire divisible par tout entier premier \`a $p$.
\end{theo}

\begin{theo} 
Soit $X$ une vari\'et\'e projective, lisse, g\'eom\'etriquement connexe sur un corps $p$-adique $k$.

{\rm (i)} Pour presque tout premier $l$, le quotient $A_{0}(X)/l$ est nul.

{\rm (ii)} Pour tout premier $l \neq p$, le quotient $A_{0}(X)/l$ est fini.
\footnote{Dans \cite{SS2}, l'\'enonc\'e (ii)  est \'etabli 
pour les $k$-vari\'et\'es qui admettent un mod\`ele   quasisemistable 
sur l'anneau des entiers; comme on verra,
le cas g\'en\'eral s'y ram\`ene gr\^ace \`a un th\'eor\`eme r\'ecent de Gabber.}
\end{theo}

Les r\'esultats de Saito et Sato \cite{SS2} furent ensuite combin\'es par Asakura et Saito \cite{AS}
\`a des techniques de th\'eorie de Hodge pour \'etablir l'existence de
surfaces $X$ de degr\'e au moins 5 dans $\P^3_{k}$ dont les sous-groupes
de torsion $l$-primaire ($l \neq p $) sont infinis (voir le  paragraphe 4 ci-apr\`es).

  \bigskip
  
  Je remercie Tam\'as Szamuely pour de nombreuses discussions  sur le th\'eor\`eme
  de Saito et Sato et pour ses commentaires critiques sur une premi\`ere version
  du pr\'esent texte.
  
  \bigskip
  
  {\bf Notations}
  
  Soit $A$ un groupe ab\'elien. Pour $n>0$ un entier, on note $A[n] \subset A$ le sous-groupe
  form\'e des \'el\'ements annul\'es par $n$. Pour $l$ un nombre premier, on note $A\{l\} \subset A$
  le sous-groupe de torsion $l$-primaire.

\section{Courbes sur les corps $p$-adiques : r\'esultats classiques}

Comme mentionn\'e dans l'introduction, pour une courbe
projective et lisse $X$ sur un corps $p$-adique, on peut d\'etecter les classes dans
$CH_{0}(X) \simeq \pic(X)$ au moyen de la cohomologie \'etale.
Expliquons plus pr\'ecis\'ement ce que nous entendons par l\`a.

\begin{theo}[Tate 1958  \cite{T}]
Soient $k$ un corps $p$-adique, 
$A$ une vari\'et\'e ab\'elienne sur $k$ et
$\hat{A}$ la vari\'et\'e ab\'elienne duale. Il y a une dualit\'e parfaite
$$ A(k) \times  H^1(k,\hat{A}) \to \br k=  \Q/\Z$$
entre le groupe ab\'elien compact $A(k)$ des points rationnels de $A$
et le groupe discret d\'efini par le premier groupe de cohomologie galoisienne
de $k$ \`a valeurs dans le groupe des points de $\hat{A}$.
\end{theo}

En s'appuyant sur ce th\'eor\`eme, on montre :

\begin{theo}[Roquette \ 1966,  \ Lichtenbaum  \ 1969  \  \cite{Li}]
\label{RoquetteLichtenbaum}
Soient $k$ un corps \mbox{$p$-adique} et $X$ une $k$-courbe projective, lisse,
g\'eom\'etriquement connexe.

{\rm a)}  Il y a un accouplement naturel 
$$ \pic X \times \br X  \to \br k=  \Q/\Z,$$
et cet accouplement est non d\'eg\'en\'er\'e des deux c\^ot\'es.

{\rm b)} Le noyau de la fl\`eche  $\Q/\Z=\br k \to \br X$
induite par le morphisme structural est $\Z/I$, o\`u $I$
est l'index de $X$, c'est-\`a-dire le pgcd des degr\'es, sur $k$,
des points ferm\'es sur~$X$.
\end{theo}

Ainsi, pour $X$ une courbe projective, lisse, g\'eom\'etriquement connexe
sur un corps~$k$ \mbox{$p$-adique,} poss\'edant un $k$-point,
le groupe $A_{0}(X)$ est isomorphe au groupe de Lie \mbox{$p$-adique} 
$J_{X}(k)$, et l'accouplement
$$ CH_{0}(X) \times \br X \to \br k = \Q/\Z$$
est non d\'eg\'en\'er\'e \`a gauche et \`a droite.

\begin{theo}[Artin 1966 \cite{G}, \S 3]
\label{ArtinGrothendieck}
Soit $R$ un anneau de valuation discr\`ete hens\'elien  excellent de corps r\'esiduel $F$.
Soit $X$ un $R$-sch\'ema fid\`element plat, projectif, int\`egre, r\'egulier,
de dimension relative 1.
Soit $X_{s}/F$ sa fibre sp\'eciale.

{\rm (a)} La fl\`eche de restriction $\br X \to \br X_{s}$ est un isomorphisme
(de groupes de torsion).

{\rm (b)} Si $F$ est un corps s\'epablement clos ou un corps  fini, $\br X \simeq \br X_{s} =0$.

{\rm (c)} Sous les m\^emes hypoth\`eses qu'en {\rm (b)}, pour tout entier $n>0$, l'injection naturelle
$$ \pic X/n \hookrightarrow H^2_{fppf}(X,\mu_{n})$$
est un isomorphisme.
\end{theo}

L'\'enonc\'e (c) provient de (b) et de la suite de Kummer
$$ 1 \to \mu_{n} \to \G_{m}  \buildrel {x \mapsto x^n} \over \longrightarrow  \G_{m} \to 1$$
sur $X$, consid\'er\'ee comme suite de faisceaux pour la topologie \'etale
sur $X$ si $n$ est inversible, et comme suite de  faisceaux pour la topologie fppf
en g\'en\'eral. La partie premi\`ere \`a $p={\rm car}(F)$ des \'enonc\'es
est plus facile \`a \'etablir, elle ne n\'ecessite pas l'hypoth\`ese d'excellence.
 
 \medskip
 
 Le th\'eor\`eme principal de Saito et Sato (th\'eor\`eme \ref{thmprincipal} ci-dessous)  g\'en\'eralise l'\'enonc\'e~(c) 
 du th\'eor\`eme \ref{ArtinGrothendieck} pour $n$ premier \`a $p$.

\section{Surfaces sur les corps $p$-adiques : quelques applications de la m\'ethode de Bloch}

Soit $X$ un sch\'ema noeth\'erien de dimension finie.
Pour tout entier
$i \geq 0$ , le groupe $Z_{i}(X)$ des cycles de dimension $i$
est le groupe ab\'elien libre sur les points (sch\'ematiques) de $X$
de dimension $i$ (ou si l'on pr\'ef\`ere  les sous-sch\'emas ferm\'es
int\`egres de dimension~$i$) 
$$Z_{i}(X)  = \oplus_{x \in X_{i}} \Z.$$
Notons $\kappa(x)$ le corps r\'esiduel en un point $x$.
On sait d\'efinir une application {\og diviseur \fg}
$$ \div :  \oplus_{x\in X_{i+1}} \kappa(x)^{\times} \to \oplus_{x \in X_{i}} \Z$$
qui g\'en\'eralise la notion de diviseur d'une fonction rationnelle (voir \cite{Ful}).
Par d\'efinition, le groupe de Chow $CH_{i}(X)$ est le conoyau de
cette application.

Supposons $X$ int\`egre et \'equidimensionnel de dimension $d$.
On note alors $CH^r(X)=CH_{d-r}(X)$.
Pour $i=d$, $CH_{0}(X)=\Z$. Pour $i=d-1$, $CH_{d-1}(X)=CH^1(X)$.
On a une fl\`eche naturelle 
$$\pic X \to CH^1(X)$$
qui est un isomorphisme si $X$ est r\'egulier.

On note $\K_{i}$  le faisceau pour la topologie de
Zariski sur $X$ associ\'e au pr\'efaisceau qui \`a un ouvert  affine $U$ associe
le groupe de $K$-th\'eorie  $K_{i}(U)$ d\'efini par Quillen.
En combinant les th\'eor\`emes sur la conjecture de Gersten,
tant en $K$-th\'eorie (Quillen) qu'en cohomologie \'etale (Bloch et Ogus \cite{BO}, 1974)
et le th\'eor\`eme de Merkur'ev et Suslin (\cite{MS}, 1982)  sur le symbole de restes normique,
 Spencer Bloch a \'etabli le r\'esultat suivant (\cite{Bl0, Bl1,Bl2}, voir aussi \cite{CT1}).

\begin{theo}[Bloch]
 \label{bloch}
Soient $k$ un corps, $X$ une $k$-vari\'et\'e lisse int\`egre,  $k(X)$ son corps des fonctions
rationnelles, et $n$ un entier
non nul dans $k$. On a une suite exacte naturelle de groupes ab\'eliens
$$0 \to H^1_{Zar}(X,\K_{2})/n \to NH^3_{\et}(X,\mu_{n}^{\otimes 2}) \to CH^2(X)[n] \to 0,$$
o\`u  $NH^3_{\et}(X,\mu_{n}^{\otimes 2})$ est le noyau de la fl\`eche de restriction
de groupes de cohomologie \'etale $H^3_{\et}(X,\mu_{n}^{\otimes 2}) \to H^3_{\et}(k(X),\mu_{n}^{\otimes 2})$.
\end{theo}

Dans la litt\'erature r\'ecente, le groupe  $H^1_{Zar}(X,\K_{2})$ a \'et\'e identifi\'e avec d'autres groupes :
le groupe de Chow sup\'erieur $CH^2(X,1)$ de Bloch d'une part,
 le   groupe d'hypercohomologie ${\bf H}^3_{Zar}(X,\Z(2))$
du complexe motivique  $\Z(2)$ d'autre part.

  Il existe des analogues de cette suite exacte
pour les sch\'emas lisses au-dessus d'un anneau de valuation discr\`ete.

Cette suite exacte  a eu de nombreuses applications,
qu'on ne saurait d\'ecrire ici  de fa\c con exhaustive.

On l'a utilis\'ee, conjointement avec le th\'eor\`eme de Deligne \'etablissant les conjectures
de Weil,
pour donner des d\'emonstrations alternatives du  th\'eor\`eme de Kato et Saito
sur le corps de classes non ramifi\'e pour les vari\'et\'es projectives et lisses
sur un corps fini mentionn\'e dans l'introduction (\cite{CTSS,Sz}).

Sur les corps $p$-adiques et sur les corps de nombres,
on l'a utilis\'ee pour obtenir des r\'esultats de finitude pour la torsion
du groupe de Chow de codimension 2, et aussi du groupe de Chow des
z\'ero-cycles, pour certaines classes de vari\'et\'es. 
Pour des r\'esultats sur les corps de nombres, je renvoie le lecteur
\`a \cite{CTR2, Sai,  Sal, PS, SS1} et aux rapports \cite{CT1, CT2}.

\begin{theo}[\cite{CTSS}]
 \label{CTSS}
Soit $k$ un corps $p$-adique.
Soit $X$ une $k$-vari\'et\'e lisse.

{\rm (i)} Pour tout entier $n>0$, le groupe $CH^2(X)[n]$
est un groupe fini.

{\rm (ii)} Pour tout $l$ premier,
le groupe de torsion $l$-primaire $CH^2(X)\{l\}$
est un groupe de cotype fini (somme d'un groupe fini $l$-primaire
et d'un groupe  $(\Q_l/\Z_l)^N$).
\end{theo}
\noindent{\sc Preuve}   --- 
D'apr\`es le th\'eor\`eme  \ref{bloch}, 
  le groupe $CH^2(X)[n]$ est un sous-quotient
du groupe $H^3_{\et}(X,\mu_n^{\otimes 2})$. La finitude de
ce groupe  pour un corps local est bien connue.
Elle implique 
que le groupe $H^3_{\et}(X,\Q_l/\Z_l(2))$ 
est un groupe de cotype fini, et donc aussi tout
sous-quotient. \cqfd

\bigskip

\begin{theo}\label{elementaire}
Soient $k$ un corps $p$-adique et $X$ une $k$-vari\'et\'e projective,  lisse, g\'eom\'etriquement int\`egre.

{\rm (i)} Pour tout $l$ premier, $l\neq p$, et tout entier $n>0$, l'application
naturelle
 $A_{0}(X)\{l\}/l^n \to A_{0}(X)/l^n$ est un isomorphisme.

{\rm (ii)} Supposons que $X/k$ a bonne r\'eduction $Y/\F$. Alors pour
tout premier $l \neq p$, l'application de r\'eduction induit une
surjection $A_{0}(X)\{l\} \to A_{0}(Y)\{l\}$.
\end{theo}

\noindent{\sc Preuve}   --- 
 (i)  
Que l'application
$A_0(X)\{l\}/l^n \to A_0(X)/l^n$  soit une injection est clair.  Pour \'etablir que c'est une
surjection,
on se ram\`ene par le th\'eor\`eme de Bertini au cas o\`u $X$ est une
courbe projective lisse g\'eom\'etriquement int\`egre, et l'assertion
r\'esulte alors de la structure du groupe des points d'une
vari\'et\'e ab\'elienne sur un corps $p$-adique.

(ii)
 On dispose d'une application de sp\'ecialisation
$ A_0(X)\to A_0(Y)$
qui est surjective (lemme de Hensel).
Soit $m>0$. D'apr\`es (i)  
 l'application $A_0(X)\{l\}/l^m \to
A_0(X)/l^m$ est surjective. L'application
compos\'ee $$A_0(X)\{l\}/l^m \to
A_0(X)/l^m \to A_0(Y)/l^m$$ est donc surjective.
Le groupe $A_0(Y)$ est fini (th\'eor\`eme de Kato et Saito \cite{KS}).
Prenant alors $m$ tel que 
 $l^m$ annule la partie $l$-primaire de $A_0(Y)$, on obtient l'\'enonc\'e.
\cqfd

\begin{theo}\label{saitosujatha}
 Soient $k$ un corps $p$-adique et
 $X$ une $k$-surface projective, lisse et
 g\'eom\'etriquement int\`egre.

{\rm (i)} Pour tout entier positif $n$, le groupe $A_0(X)[n]$
est un groupe fini.

{\rm (ii)} Pour tout $l$ premier,
le groupe de torsion $l$-primaire $A_0(X)\{l\}$
est un groupe de cotype fini (somme d'un groupe fini $l$-primaire
et d'un groupe  $(\Q_l/\Z_l)^N$).

{\rm (iii) (Saito et Sujatha)} Pour $l$ premier, $l\neq p$, le groupe $A_0(X)$
est somme de son sous-groupe $l$-divisible maximal et d'un groupe
fini $l$-primaire.

{\rm (iv) (Saito et Sujatha)} Pour $n>0$ premier \`a $p$, le groupe $A_0(X)/n$ est fini.
\end{theo}
\noindent{\sc Preuve} (voir \cite{CT2})  --- 
 Pour $X$ une surface projective lisse g\'eom\'etriquement
connexe, on a l'\'egalit\'e $CH_0(X)=CH^2(X)$.
Les \'enonc\'es (i) et (ii) sont
des cas particuliers du th\'eor\`eme \ref{CTSS}.
D'apr\`es (ii) on peut
\'ecrire
$A_0(X)\{l\}=(\Q_l/\Z_l)^N \oplus F_l$ avec
$F_l$ un groupe ab\'elien fini annul\'e par une puissance $l^t$ de $l$,
et $N\geq 0$. En utilisant la proposition \ref{elementaire}(i) on voit alors
que
l'application compos\'ee $F_l \to A_0(X) \to
A_0(X)/l^t$ est  un isomorphisme. On en d\'eduit que dans la suite
exacte
$$0 \to D_l \to A_0(X) \to A_0(X)/l^t \to 0$$
d\'efinissant $D_l$, la projection $A_0(X) \to A_0(X)/l^t$
est scind\'ee, d'o\`u $A_0(X) \simeq D_l \oplus F_l,$
avec $D_l/l=0$. Le groupe $D_l$ est donc le sous-groupe
$l$-divisible maximal de $A_0(X)$. Ceci \'etablit (iii), et 
l'\'enonc\'e (iv) suit. \cqfd

\begin{rema}
On verra au paragraphe 3 que pour presque tout premier $l$ le groupe fini $F_{l}$ est nul,
en d'autres termes le groupe $A_{0}(X)$ est $l$-divisible.
La m\'ethode ci-dessus ne permet pas d'obtenir ce r\'esultat.
\end{rema}

\bigskip

Le th\'eor\`eme suivant, d\'etaill\'e dans \cite{CT3}, regroupe des travaux des ann\'ees 1985 \`a 1991, dus \`a 
Raskind et au r\'edacteur \cite{CTR1,CTR2},  \`a Salberger \cite{Sal} et \`a S. Saito \cite{Sai}.

\begin{theo}
 Soient $k$ un corps $p$-adique et
$X$ une
$k$-surface projective et lisse, g\'eom\'etriquement int\`egre.
 Supposons $H^2(X,O_X)=0$.  Alors

{\rm (i)} Le groupe $A_0(X)_{tors}$ est fini.

{\rm (ii)}  Si la fl\`eche naturelle $A_0({\overline X}) \to
\Alb_X({\overline k})$ (sur une cl\^oture alg\'ebrique $\overline k$ de $k$)
est un isomorphisme (ce qui r\'esulterait
de $H^2(X,O_X)=0$ suivant une conjecture de S.~Bloch, connue pour les surfaces qui ne sont pas de type g\'en\'eral), alors
le groupe $A_0(X)$ est extension d'un sous-groupe ouvert
de $\Alb_X(k)$ par un groupe fini. En particulier, 
le groupe $A_0(X)$ est somme directe d'un groupe fini d'ordre premier \`a $p$
et d'un groupe $p'$-divisible.

{\rm (iii)}  Si la  vari\'et\'e d'Albanese de $X$  a bonne r\'eduction, 
alors l'accouplement $$A_0(X)_{tors} \times \br(X) \to \Q/\Z $$ 
est non d\'eg\'en\'er\'e \`a gauche.

{\rm (iv)} Si la fl\`eche naturelle $A_0({\overline X}) \to
\Alb_X({\overline k})$ 
est un isomorphisme et si la  vari\'et\'e d'Albanese de $X$  a bonne r\'eduction, 
alors l'accouplement $$A_0(X)
\times \br(X) \to \Q/\Z $$ est non d\'eg\'en\'er\'e \`a gauche. 
\end{theo}

L'\'enonc\'e (iii) peut \^etre \'etabli sous des hypoth\`eses
un peu plus larges (\cite{Sal}, \cite{Sai}, \cite{Sato1}, \cite{SS1}), mais
on ne peut  totalement ignorer les  hypoth\`eses   dans (iii) et (iv).
Parimala et Suresh \cite{PS} construisent une surface $X$
fibr\'ee en coniques lisses au-dessus d'une
courbe $C$ (la surface satisfait donc $H^2(X,O_{X})=0$),
la conique ayant mauvaise r\'eduction sur un corps $p$-adique $k$
de caract\'eristique  r\'esiduelle impaire, surface
 pour laquelle :

(a) le noyau \`a gauche de l'accouplement
$A_0(X)\{2\}   \times \br X \to \Q/\Z $ 
n'est pas nul;

(b) l'application cycle $CH^2(X)/2 \to H^4_{\et}(X,\Z/2)$ n'est pas injective;

(c) le noyau de l'application
$ A_{0}(X) \to \Hom(\br X, \Q/\Z)$ n'est pas le sous-groupe divisible maximal
de $A_{0}(X)$.

\bigskip

La plupart des th\'eor\`emes pr\'ec\'edents furent obtenus en \'etudiant l'action du groupe
de Galois absolu de $k$ sur divers groupes ($\K$-cohomologie, cohomologie \'etale) 
 attach\'es aux vari\'et\'es (apr\`es passage
\`a une cl\^oture alg\'ebrique du corps $k$).

Une autre m\'ethode consiste \`a consid\'erer des mod\`eles des $k$-vari\'et\'es
au-dessus d'un ouvert de l'anneau des entiers de $k$. Limitons-nous ici
au  cas de bonne r\'eduction.
Soit $R$ l'anneau des entiers d'un corps $p$-adique
$k$, de corps r\'esiduel $\F$. Soit  $\cal X$ un $R$-sch\'ema int\`egre,
projectif et lisse, de fibre g\'en\'erique $X/k$ g\'eom\'etriquement
int\`egre, et soit $Y/\F$ la fibre sp\'eciale.
Pour un tel $\X$,  on a la suite exacte de localisation
$$ H^1(X,{\cal K}_2) \to \pic(Y) \to CH^2({\cal X}) \to CH^2(X) \to 0.$$
Notons $\delta : H^1(X,{\cal K}_2) \to \pic(Y)$.
Le groupe $\pic(Y)$ est un groupe de type fini (th\'eor\`eme de
N\'eron-Severi et finitude du groupe des points rationnels
d'une vari\'et\'e ab\'elienne sur un corps fini). 
L'\'enonc\'e suivant, d\'etaill\'e dans \cite{CT3},  regroupe des r\'esultats de 
Raskind \cite{R}, Raskind et l'auteur \cite{CTR2}, Spie{\ss} \cite{Sp}.

\begin{theo}
Avec les notations ci-dessus, supposons $\dim(X)=2$, et
supposons

{\rm (H)} L'application $H^1(X,{\cal K}_2)\otimes \Q  \to \pic(Y)\otimes \Q$ est surjective.

 Alors

{\rm (i)} L'application de sp\'ecialisation 
$CH_0(X) \to CH_0(Y)$
induit un isomorphisme sur les sous-groupes de torsion
premi\`ere \`a $p$.

{\rm (ii)} Le groupe $A_0(X)$ est la somme directe d'un groupe
fini d'ordre premier \`a $p$ et d'un groupe uniquement divisible
par tout entier $n$ premier \`a $p$. Pour tout $l$ premier
distinct de $p$, le groupe $D_l$ ci-dessus est uniquement
$l$-divisible.

{\rm (iii)} Pour $n>0$ premier \`a $p$, l'application cycle
$$CH^2(X)/n \to H^4_{\et}(X,\mu_n^{\otimes 2})$$
est injective.

{\rm (iv)} L'accouplement
$$ A_0(X) \times \br(X) \to \Q/\Z$$
a son noyau \`a gauche form\'e d'\'el\'ements divisibles
par tout entier $n>0$ premier \`a $p$.
\end{theo}

Voici des cas o\`u l'hypoth\`ese (H) a \'et\'e \'etablie.

a) Le cas o\`u $H^2(Y,O_Y)=0$ (et donc aussi $H^2(X,O_X)=0$).
C'est le cas le plus simple.
Dans ce cas le conoyau de la fl\`eche compos\'ee 
$\pic(X)\otimes k^{\times} \to H^1(X,{\cal K}_2) \to \pic(Y)$
est nul, car la fl\`eche de restriction $\pic({\cal X}) \to \pic(Y)$ est
surjective. Ce cas fut consid\'er\'e par Raskind \cite{R}, Coombes,
CT-Raskind (\cite{CTR2}).

b) Le cas o\`u  le rang du groupe
de N\'eron-Severi g\'eom\'etrique ne grandit pas par sp\'ecialisation
(Raskind \cite{R}).

c) Le cas o\`u ${\cal X}$ est le produit fibr\'e de deux
courbes elliptiques avec bonne r\'eduction (Spie{\ss} \cite{Sp}).

\bigskip

Dans le cas c) il faut, pour \'etablir (H), trouver des
\'el\'ements  \og ind\'ecomposables \fg  \  dans $H^1(X,{\cal K}_2 )$, i.e.
d'autres
\'el\'ements   que ceux  provenant de $\pic(X) \otimes k^{\times}$.
Spie{\ss} utilise certaines correspondances entre courbes
elliptiques provenant de travaux de Frey et Kani.

Dans les ann\'ees 1990 \`a 2000, il y eut dans cette direction (recherche d'\'el\'ements ind\'ecomposables) 
une s\'erie de travaux
(Flach, Mildenhall,   S. Saito, Langer, Raskind, Otsubo)
reposant souvent sur une arithm\'etique
tr\`es fine des vari\'et\'es consid\'er\'ees (par exemple des produits
de courbes modulaires).  Je renvoie ici le
 lecteur \`a l'article r\'ecent  \cite{SS1} de S. Saito et K. Sato, tant pour la liste
 de r\'ef\'erences que pour  le lien entre la finitude
 du groupe $CH^2(X)_{tors}$ et le comportement 
 d'applications
 \og r\'egulateurs \fg \  de source  le groupe  $H^1(X,\K_{2})$, applications
 d\'ej\`a consid\'er\'ees par Salberger   \cite{Sal}.

Le r\'ecent th\'eor\`eme de Asakura et Saito (Th\'eor\`eme \ref{asakurasaito} ci-dessous)
montre que
 l'hypoth\`ese (H) ne vaut pas toujours : elle est en d\'efaut pour des
 surfaces lisses dans $\P^3_{k}$ de degr\'e au moins 5, suffisamment g\'en\'eriques.
Asakura et Sato 
ont pos\'e la question de sa validit\'e 
 lorsque la surface sur le corps $p$-adique  provient
 d'une surface d\'efinie sur un corps de nombres.

\section{LE TH\'EOR\`EME DE S. SAITO ET K. SATO  \cite{SS2}}

Dans tout ce paragraphe, on adopte les notations suivantes.

On note $R$ un anneau de valuation discr\`ete,  
$F$ son corps
r\'esiduel et $k$ son corps des fractions. On note $B=\spec R$, $s=\spec F$, 
$\eta = \spec k$. 

On note $Sch^{qp}_{B}$ la cat\'egorie des sch\'emas quasi-projectifs sur $B$.
Pour $X \in Sch^{qp}_{B}$, on note $\delta(X) \in \N$ la dimension de Krull d'une compactification
de $X$ au-dessus de $B$.

Pour $X \in Ob(Sch^{qp}_{B})$, on note $X_{s}/F$ sa fibre sp\'eciale et $X_{\eta}/k$ sa fibre g\'en\'erique.
Si $X$ est int\`egre et $X_{\eta}$ est vide, donc $X=X_{s}$,
alors $\delta (X)$ est \'egal au degr\'e de transcendance sur $F$ du corps des fonctions 
 $F(X_{s})$.  Si $X$ est int\`egre et $X_{\eta}$ est non vide, alors $\delta (X)-1$
 est \'egal au degr\'e de transcendance du corps des fonctions $k(X_{\eta})$ sur $k$.

On note $\cal C$ la sous-cat\'egorie  pleine de  $Sch^{qp}_{B}$ dont les objets
satisfont $X_{s} \neq \emptyset$.
Pour $X$ irr\'eductible dans $\mathcal C$, on a
$\delta (X)= \dim(X)$.

Un objet $X \in \mathcal{C}$ est appel\'e quasi-semistable s'il satisfait les conditions :

(QS1) $X$ est r\'egulier, \'equidimensionnel,  plat et de type fini sur $B$.

(QS2) Le diviseur r\'eduit $X_{s,red}$ sur $X$ est \`a croisements normaux stricts.

On note $\mathcal{QS} \subset \mathcal{C}$, resp. $\mathcal{QSP} \subset \mathcal{C}$,
la sous-cat\'egorie pleine dont les objets sont les objets quasi-semistables,
resp. les objets quasi-semistables et projectifs sur $B$.

Une $\mathcal{QS}$-paire est un couple ($X,Y)$ de sch\'emas dans $Ob(\mathcal{QS})$
pour lequel $Y$ est un diviseur sur $X$ et le diviseur $X_{s,red} \cup Y$
sur $X$ est \`a croisements normaux stricts. 
Posant $U: = X \setminus Y$, on  note indiff\'eremment $(X,Y)=(X,Y;U)$.

Une $\mathcal{QSP}$-paire est une $\mathcal{QS}$-paire   $(X,Y)$ pour laquelle $X$ et $Y$ sont
projectifs sur $B$.

Une $\mathcal{QSP}$-paire ample   est une  $\mathcal{QSP}$-paire $(X,Y;U)$
pour laquelle $U$ est affine.

\subsection{Groupes de Chow, homologie \'etale, suite spectrale de niveau,
complexe de Bloch-Ogus et Kato,
 application cycle}

Pour les vari\'et\'es alg\'ebriques  sur un corps,
 les r\'esultats expos\'es au paragraphe 2 utilisent diverses propri\'et\'es
 de la cohomologie et de  l'homologie \'etale, la dualit\'e de Poincar\'e,
 les filtrations par le niveau
 et par le coniveau et les suites spectrales associ\'ees.
 Ces th\'eories sont d\'evelopp\'ees dans 
 \cite{G} III, \S 10.1, \cite{La}, \cite{BO}.
 
 Il a fallu \'etendre ces th\'eories
 aux sch\'emas de type fini au-dessus d'un anneau de valuation discr\`ete.
 
 Sauf mention du contraire, la cohomologie employ\'ee est la cohomologie \'etale.

\medskip

Pour $X$ dans $\mathcal C$, et $q\geq 0$ entier, on note $X_{q}$ l'ensemble des points
 $x \in X$ dont l'adh\'erence $\overline{\{x\}} \subset X$ satisfait $\delta(\overline{\{x\}} ) =q$.  
Pour $q \geq 0$, et $X \in  {\cal C}$, on d\'efinit le groupe de Chow de dimension $q$ par
la formule usuelle :
$$CH_{q}(X)= {\rm Coker} \,  [\div : \oplus_{x \in X_{q+1}}  \kappa(x)^* \to \oplus_{x \in X_{q}} \Z].$$

On fixe un premier $l \neq {\rm car} \ F$. 
On note $\mu_{l^n}$ le faisceau pour la topologie \'etale associ\'e au  $B$-sch\'ema 
en groupes des racines $l^{n}$-i\`emes de l'unit\'e.
Pour $m \in \N$, on a le faisceau $\Z/l^n(m)= \mu_{l^n}^{\otimes m}$.
Pour $-m \in \N$ on note  $\Z/l^n(m)$ le faisceau $\Hom(\Z/l^n(-m),\Z/l^n).$
Pour $m \in \Z$ on note $\Q_{l}/\Z_{l}(m)=\colime\Z/l(m)$.
Lorsque l'on voudra consid\'erer simultan\'ement le cas $\Z/l^n$ et le cas $\Q_{l}/\Z_{l}$,
on utilisera la notation $\Lambda$.

 Pour $X \in  Ob(Sch^{qp}_{B})$
de morphisme structural $f : X \to B$ 
et pour $q \in \Z$ et $m \in \Z$ ,  suivant Grothendieck, Artin, Verdier, Deligne,  
on d\'efinit  
$$H_{q}(X,\Lambda(m)) : = {\mathbb H}^{2-q}(X,Rf^{!}\Lambda(m))$$
o\`u $Rf^{!}$ est le foncteur image inverse extraordinaire ([SGA4, XVIII, Thm. 3.1.4]).

Ceci d\'efinit une
  th\'eorie homologique sur la cat\'egorie $Sch^{qp}_{B}$
  ayant toutes les propri\'et\'es voulues :  fonctorialit\'e covariante 
  par morphisme propre, fonctorialit\'e contravariante par morphisme \'etale,
  existence pour toute immersion ferm\'ee $Y \hookrightarrow X$ de compl\'ementaire
  $U \hookrightarrow X$  d'une longue suite exacte
  $$ \dots \to H_{q}(Y,\Lambda) \to H_{q}(X,\Lambda)  \to H_{q}(U,\Lambda)  \to H_{q-1}(Y,\Lambda) \to \dots,$$
  fonctorialit\'e de cette suite exacte. Le lecteur se reportera \`a \cite{La, JS1, JSS}.

  Cette th\'eorie homologique \`a la Borel-Moore se relie \`a la cohomologie \'etale gr\^ace
  \`a un th\'eor\`eme de \og dualit\'e de Poincar\'e \fg :
  
  \begin{prop}\label{dualite}
  Soit  $X \in  Ob(Sch^{qp}_{B})$ un sch\'ema r\'egulier int\`egre. Posons $d=\delta (X)-1$.
  Pour tout sous-sch\'ema ferm\'e $Y \subset X$, et tout $q \in \Z$,
   il existe un isomorphisme canonique
  $$ H^{2d-q+2}_{Y}(X,\Lambda(d)) 
  \xrightarrow{\cong}
  H_{q}(Y,\Lambda).$$
    \end{prop}
  \noindent    (Par convention $H^r=0$ pour $r<0$.)

   Cet isomorphisme satisfait une s\'erie de propri\'et\'es fonctorielles (cf.  \cite{SS2}).
     
Pour $X \in  Ob(Sch^{qp}_{B})$, $x$ un point de $X$,  et $q \in \Z$   on note
$$ H_{q}(x,\Lambda): = \colim_{U \subset {\overline{\{x\}}   }} H_{q}(U,\Lambda),$$
o\`u ${\overline{\{x\}}}$ est l'adh\'erence de $x$ dans $X$ et $U$ parcourt les ouverts non vides de ${\overline{\{x\}}}$.

La longue suite exacte de localisation donn\'ee ci-dessus et la filtration par le niveau
donnent naissance \`a une suite spectrale de type homologique
$$E^1_{a,b}(X,\Lambda) =  \oplus_{x\in X_{a}} H_{a+b}(x,\Lambda) \Longrightarrow  
H_{a+b}(X,\Lambda)$$
 dont  les diff\'erentielles de niveau $r$ sont de degr\'e $(-r,r-1)$.

Pour $n \in \Z$ et $m \in \Z$, notons $H^n(x,\Lambda(m))$ le groupe de cohomologie galoisienne 
$H^n(\kappa(x),\Lambda(m))$, groupe qui par d\'efinition est nul pour $n<0$.
Par passage \`a la limite dans les isomorphismes dans la proposition \ref{dualite} 
(pour $Y=X=U$)
on obtient la

\begin{prop}[Suite spectrale de niveau]\label{ssniveau}
Pour $X \in  Ob(Sch^{qp}_{B})$, il y a une suite spectrale homologique
$$E^1_{a,b}(X,\Lambda) =  \oplus_{x\in X_{a}} H^{a-b}(x,\Lambda(a-1)) \Longrightarrow  H_{a+b}(X,\Lambda).$$
Cette suite exacte est fonctorielle covariante par rapport aux morphismes propres et contravariante par
rapport aux morphismes \'etales.
\end{prop}

La ligne $b=0$ de la suite spectrale est un complexe 
$$0 \leftarrow \bigoplus_{x \in X_{0}} H^0(x,\Lambda(-1)) \leftarrow \bigoplus_{x \in X_{1}} H^1(x,\Lambda)
\leftarrow \cdots \bigoplus_{x \in X_{a}} H^a(x,\Lambda(a-1)) \leftarrow  \cdots $$
(la somme sur les points de dimension $a$ \'etant plac\'ee en degr\'e $a$). 
Pour $X \in \cal{C} $ int\`egre, de corps des fonctions $\kappa(X)$, 
avec $d=\delta (X)-1={\rm dim}(X)-1$, le complexe commence (\`a droite)  par
$$  \cdots \leftarrow H^{d+1}(\kappa(X),\Lambda(d))  \leftarrow 0.$$
On note ce complexe $KC(X,\Lambda)$. On note $KH_{a}(X,\Lambda)$ le groupe d'homologie en degr\'e~$a$.
D'apr\`es Jannsen, Saito et Sato \cite{JSS}, les fl\`eches dans ce complexe sont les oppos\'ees
des fl\`eches de bord en cohomologie galoisienne utilis\'ees par Kato \cite{Kato}.

\begin{rema}
Lorsque l'on \'etudie  ([K], \cite{J}, \cite{JS3}) le corps de classes
de dimension sup\'erieure sur le corps des fonctions d'une vari\'et\'e int\`egre $X$ 
de dimension $d+1$
sur un corps fini  $\F$, le principal complexe consid\'er\'e  va de $H^{d+2}(\kappa(X),\Lambda(d+1)) $ \`a
$\bigoplus_{x \in X_{0}} H^1(x,\Lambda)$.
\end{rema}

Dans la suite spectrale de niveau, on a clairement $E^1_{a,b}(X,\Lambda)=0$ pour $a \notin [0, \delta(X)]$ et pour $a-b <0$.

On a en particulier des applications 
$$E^1_{1,1} \to E^2_{1,1} \to E^{\infty}_{1,1} \hookrightarrow H_{2}(X,\Lambda).$$

Le groupe $E^2_{1,1}$ est le conoyau de la fl\`eche $d^1_{2,1} : E^1_{2,1} \to E^1_{1,1}$,
donc (en utilisant la suite de Kummer)
$$E^2_{1,1}(X,\Lambda) = {\rm Coker} \,  [d^1_{2,1} :   \oplus_{x \in X_{2}} \kappa(x)^{\times} \otimes \Lambda
\to \oplus_{x \in X_{1}} \Lambda].$$
Les auteurs identifient cette application avec l'oppos\'ee de l'application 
 $\div \otimes \Lambda$.
Pour $X \in \cal C$, le conoyau s'identifie donc avec le
 groupe $ CH_{1}(X)  \otimes \Lambda$. 
L'application compos\'ee $ E^2_{1,1} \to E^{\infty}_{1,1} \hookrightarrow H_{2}(X,\Lambda)$ 
d\'efinit un homomorphisme
$$ \rho_{X } : CH_{1}(X)  \otimes \Lambda \to H_{2}(X,\Lambda)$$
dont on v\'erifie 
qu'elle co\"{\i}ncide
avec l'application cycle (sur les $1$-cycles de $X$).  
Cette application jouit des propri\'et\'es fonctorielles attendues :
elle est covariante par morphismes propres dans $\cal C$, et contravariante par morphismes \'etales
dans $\cal C$. 
 
 La fonctorialit\'e suivante est particuli\`erement importante.
Pour une immersion ferm\'ee $Y \subset X$ dans $\cal C$ dont le compl\'ementaire
$U \hookrightarrow X$ est dans $\cal C$, on a un diagramme commutatif 
de suites de localisation
  \[
\xymatrix{
CH_{1}(Y)\otimes \Lambda \ar[r]  \ar[d]_{\rho_{Y}} &  CH_{1}(X)\otimes \Lambda   \ar[r] 
  \ar[d]_{\rho_{X}}  & CH_{1}(U)\otimes \Lambda  \ar[d]_{\rho_{U}}  \ar[r] & 0
 \\
H_{2}(Y,\Lambda) \ar[r]  & H_{2}(X,\Lambda)  \ar[r]  & H_{2}(U,\Lambda ) &
}
\]
dont les lignes sont exactes.

\begin{prop}\label{immediat}
Soit
$X \in Ob(\cal{C})$.

{\rm (1)} Si   $\delta (X)=1$,  la fl\`eche
$\rho_{X} : CH_{1}(X)\otimes \Lambda \to H_{2}(X,\Lambda)$
est un isomorphisme.

{\rm (2)}  Si $\delta (X)=2$, on a une suite exacte
$$ 0 \to CH_{1}(X)\otimes \Lambda \to H_{2}(X,\Lambda) \to KH_{2}(X,\Lambda) \to 0.$$

{\rm (3)} Si $\delta (X)=3$, on a une suite exacte
$$ H_{3}(X,\Lambda) \to KH_{3}(X,\Lambda) \to CH_{1}(X)\otimes \Lambda \to H_{2}(X,\Lambda).$$

{\rm (4)} Si $\delta (X) \leq 3$, on  a $KH_{3}(X,\Z/l^n) = KH_{3}(X,\Q_{l}/\Z_{l}) [l^n]$.
\end{prop}
\noindent{\sc Preuve} --- Les \'enonc\'es (1) \`a (3) r\'esultent imm\'ediatement de la forme de la suite
spectrale de niveau. L'\'enonc\'e (4)  r\'esulte  de la fonctorialit\'e en les coefficients de la suite spectrale, 
du th\'eor\`eme 90 de Hilbert, et du th\'eor\`eme de Merkur'ev--Suslin \cite{MS}.
\cqfd

\begin{theo}\label{dimension2}
Supposons  $R$ hens\'elien et $F$ fini ou s\'eparablement clos. 
Soit $X$ un $R$-sch\'ema r\'egulier, projectif et plat sur $R$,
de dimension 2.
Soit 
 $\Lambda=\Z/l^n $ ou $\Lambda=\Q_{l}/\Z_{l}$.
 
 Alors
 
{\rm  (a)}  L'application cycle $\rho_{X} : CH_{1}(X) \otimes \Lambda \to H_{2}(X,\Lambda)$
est bijective.

{\rm (b)} Le groupe d'homologie $KH_{2}(X,\Lambda)=0$.
\end{theo}
\noindent{\sc Preuve}   ---  On peut supposer $X$ int\`egre et $\Lambda=\Z/l^n$.
Pour  $\dim(X)=2$, la suite exacte
$$ 0 \to CH_{1}(X)/l^n  \to H_{2}(X,\Z/l^n) \to KH_{2}(X,\Z/l^n) \to 0$$
de la proposition \ref{immediat}  s'identifie \`a  
la suite  exacte 
$$ 0 \to \pic(X)/l^n \to H^2(X,\mu_{l^n}) \to \br(X)[l^n] \to 0$$
 d\'eduite de la suite de Kummer en cohomologie \'etale.
 L'\'enonc\'e est alors le th\'eor\`eme d'Artin (Th\'eor\`eme \ref{ArtinGrothendieck}).
 \cqfd

\subsection{Deux conjectures}

Soient  $X \to B$ dans  $ \cal C$ et $x \in X_{a}$. Si $p(x)=s$, alors
$cd_{l}(\kappa(x))=cd_{l}(F)+a$. Si $p(x)=\eta$, alors
$cd_{l}(\kappa(x))=cd_{l}(k)+a-1$, ce qui si $R$ est hens\'elien
implique $cd_{l}(\kappa(x))=cd_{l}(F)+a$.

Ceci implique imm\'ediatement
 les \'enonc\'es  (i) et (ii) dans la proposition
suivante. Cette proposition n'est pas utilis\'ee 
dans la d\'emonstration
du th\'eor\`eme principal \ref{thmprincipal}, elle ne sert que dans la d\'emonstration du
th\'eor\`eme \ref{annulationKH3}, qui \'etablit les conjectures \ref{conjsepclos}  et \ref{conjfini}
en bas degr\'e.

\begin{prop}\label{annulationgaloisienne}
Supposons 
$R$ hens\'elien. Soit $X \in Ob({\cal C})$.  

{\rm (i)} Si $F$ est s\'eparablement clos,  alors $E^1_{a,b}(X,\Z/l^n) =0$ pour $b <0$ et $a$ quelconque.

{\rm (ii)} Si $F$ est un corps fini, alors $E^1_{a,b}(X,\Z/l^n) =0$ pour $b <-1$ et $a$ quelconque.

{\rm (iii)} Si $F$ est un corps fini, alors $E^1_{a,-1}(X,\Q_{l}/\Z_{l}) =0$ pour tout $a$.
\end{prop}

L'\'enonc\'e (iii), qui est un cas particulier d'un \'enonc\'e d'annulation pour la cohomologie
galoisienne \`a coefficients $\Q_{l}/\Z_{l}(m)$ pour certaines torsions $m$,
 est moins classique.  Un tel \'enonc\'e avait \'et\'e obtenu par B. Kahn  \cite{Kahn}.   
Dans les situations de (i)  et (iii),   
les termes non nuls de la suite spectrale sont tous dans
le premier quadrant.
Dans le cas~(ii), il  y a une ligne suppl\'ementaire en dessous du premier quadrant.

Supposons $R$ hens\'elien excellent et $F$ s\'eparablement clos.
Pour $X \in \cal{C}$, r\'egulier, projectif,  le complexe $KC(X,\Z/l^n)$
est exact en degr\'e $a=0, 1$. Ceci r\'esulte de la forme de la suite spectrale, de l'isomorphisme
$H^n(X,\Z/l^n) \simeq H^n(X_{s},\Z/l^n)$ (th\'eor\`eme de changement de base)
et de la nullit\'e de $H^n(X_{s},\Z/l^n)$ pour $n>2d=2 \dim(X_{s})$.

Inspir\'es par des conjectures de Kato \cite{Kato} sur les vari\'et\'es sur les corps finis,
Saito et Sato \cite{SS2} sugg\`erent :

\begin{conj}\label{conjsepclos}
 Supposons $R$ hens\'elien excellent et $F$ s\'eparablement clos.
Pour tout $X \in Ob({\cal QSP})$
 le complexe $KC(X,\Z/l^n)$ est exact.
 \end{conj}

\begin{conj}\label{conjfini}
 Supposons $R$ hens\'elien excellent et $F$ fini.
Notons $I(X_{s})$   l'ensemble des composantes irr\'eductibles de $X_{s,red}$.
Pour tout $X \in Ob({\cal QSP})$
 le complexe $KC(X,\Q_{l}/\Z_{l})$ est exact sauf en degr\'e $a=1$,
 o\`u l'homologie $KH_{1}(X,\Q_{l}/\Z_{l})$    s'identifie naturellement \`a $\Q_{l}/\Z_{l}^{I(X_{s})}$.
  \end{conj}
 
La forme de la suite spectrale et des r\'esultats classiques de cohomologie \'etale [SGA4] montrent ici encore que
  cette conjecture vaut en degr\'e $a=0,1$.

\subsection{Th\'eor\`eme de Lefschetz affine affin\'e}

\begin{theo}\label{Lefschetzaffine}
Soit $(X,Y;U)$ une  $\cal{QSP}$-paire ample avec  $\dim (X) = d+1 \geq 2$.
Supposons $R$  hens\'elien.

{\rm (i)} Si $F$ est s\'eparablement clos, alors $H_{q}(U,\Z/l^n)=0$ pour
$q \leq d+1$.

{\rm (ii)} Si $F$ est fini, alors  $H_{q}(U,\Z/l^n)=0$ pour $q \leq d$.

{\rm (iii)} Si $F$ est fini et $\dim(X) \geq 3$, alors $H_{d+1}(U,\Q_{l}/\Z_{l})=0$.
\end{theo}

\noindent{\sc Preuve}  (esquisse)  ---
Le sch\'ema $U$ est r\'egulier, plat et de dimension relative $d$ sur $B$.
D'apr\`es la proposition \ref{dualite},  
pour tout entier $q$, on a un isomorphisme 
$$ H_{q}(U,\Lambda) \simeq H^{2d+2-q}(U,\Lambda(d)).$$

 Le sch\'ema $U$ est affine et $R$ est hens\'elien. Les th\'eor\`emes
 de Lefschetz affines de Artin et Gabber (voir Fujiwara \cite{Fuj}) 
 donnent la
 nullit\'e de ces groupes pour $q \leq d$ si $F$ est s\'eparablement clos
 et pour $q \leq d-1$ si $F$ est fini.  C'est un de moins que ce qui
 est affirm\'e dans (i) et (ii). Pour aller plus loin, il faut utiliser le
 fait que $(X,Y;U)$ est une $\cal{QSP}$-paire, en particulier 
 que $X_{s,red} \cup Y$ est un diviseur \`a croisements normaux stricts.

En utilisant le th\'eor\`eme de puret\'e cohomologique absolu de Gabber,
on \'etend  le th\'eor\`eme de changement de base
de Rapoport et Zink   \cite{RZ}
aux $\cal{QSP}$-paires, et l'on montre :
 \begin{prop}
Soit $(X,Y;U)$ une $\cal{QSP}$-paire sur $B$. Soit $V:=U_{s,red}$.
 Si $R$ est hens\'elien,
pour tous entiers $q$ et $j$, la fl\`eche
  image r\'eciproque $$H^q(U,\Z/l^n(j)) \to H^q(V,\Z/l^n(j))$$
  est un isomorphisme.
\end{prop}

Combinant cette proposition et 
 l'isomorphisme ci-dessus, pour \'etablir le th\'eor\`eme 
 on est ramen\'e \`a \'etablir l'annulation des groupes $H^{2d+2-q}(V,\Lambda(d))$
 pour $q$ comme dans le th\'eor\`eme.
  Les \'enonc\'es (i) et (ii) 
 sont des cons\'equences du  th\'eor\`eme
de Lefschetz affine ([SGA4], XIV, Cor. 3.2). Pour l'\'enonc\'e (iii), on montre $H^{d+1}(V,\Q_{l}/\Z_{l}(d))=0$ pour 
$F$ fini et $d \geq 2$.
En tenant compte de l'\'enonc\'e  (ii) pour $V$, on est ramen\'e \`a montrer $H^{d+1}(V,\Q_{l}(d))=0$
puis \`a voir que l'op\'erateur de Frobenius sur $H^{d+1}(V\times_{F}{\overline F},\Q_{l}(d))$
n'a pas la valeur propre 1. Ceci r\'esulte des hypoth\`eses et des th\'eor\`emes de Deligne
sur la cohomologie des vari\'et\'es sur les corps finis.

\subsection{Th\'eor\`emes de Bertini relatifs (Jannsen et Saito)}

La litt\'erature contient diverses versions du th\'eor\`eme de Bertini sur les sections hyperplanes
d'une vari\'et\'e $X$  projective et lisse sur un corps $k$.

Dans un premier temps, on affirme l'existence de sections hyperplanes $X \cap H$ lisses sur $k$.
Une r\'ef\'erence sur le sujet est un livre de J.-P.~Jouanolou.
Ceci vaut si le corps $k$  est infini.
Si le corps est un corps fini $\F$,
l'\'enonc\'e vaut encore si l'on remplace le plongement projectif donn\'e 
par un plongement de Veronese convenable.
B.~Poonen et O.~Gabber ont donn\'e des \'enonc\'es pr\'ecis dans cette direction. Pour les probl\`emes
ici consid\'er\'es, on peut aussi garder le plongement donn\'e et remplacer le corps fini $\F$  donn\'e
par deux extensions finies de $\F$ suffisamment grandes et de degr\'es premiers entre eux.

Dans un second temps, on se donne une sous-$k$-vari\'et\'e ferm\'ee  $Z \subset X$
et l'on demande s'il existe une section hyperplane lisse de $X$ qui contient $Z$.
Il faut bien s\^ur imposer des conditions \`a $Z$, par exemple
$Z$   lisse et $2\dim(Z) < \dim(X)$.
A.~Altman et S.~Kleiman ont consacr\'e un article \`a cette question.
Sur un corps fini, Poonen a obtenu les \'enonc\'es ad\'equats.

Quand on \'etudie les sch\'emas arithm\'etiques, 
pour faire des d\'emonstrations par r\'ecurrence sur la dimension, 
 on a besoin d'adapter ces th\'eor\`emes
au-dessus d'un anneau de valuation discr\`ete (ou mieux d'un sch\'ema de Dedekind).

L'\'enonc\'e suivant (Jannsen et Saito \cite{JS2}) correspond \`a la premi\`ere situation.

\begin{theo}[Th\'eor\`eme de Bertini relatif, version 1]
\label{bertini1}
Soit $R$ un anneau de valuation discr\`ete. Soit $X \in Ob(\cal{QSP})$ de dimension 
$\dim(X) \geq 2$. Il existe un \mbox{$R$-plongement} $X \hookrightarrow \P^n_{R}$
et un $R$-hyperplan $H=\P^{n-1}_{R} \subset  \P^n_{R}$ tels que l'intersection
sch\'ematique $X \cap H $ soit dans $Ob(\cal{QSP})$ et que le couple 
$(X, X\cap H)$ soit une  $\cal{QSP}$-paire (ample).
\end {theo}

La seconde situation m\`ene \`a un \'enonc\'e  plus technique, d\^u aussi
\`a Jannsen et Saito (\cite{SS2}, Thm. 4.2).  On impose en particulier des conditions de transversalit\'e
entre le sous-sch\'ema r\'egulier $Z \subset X$ consid\'er\'e et les composantes de la fibre
sp\'eciale r\'eduite.
 Je ne  reproduis pas cet \'enonc\'e ici.
Dans ce rapport, je me contente de l'appeler  \og {\it Th\'eor\`eme de Bertini relatif, version 2} \fg.

\subsection{Lemme de d\'eplacement;  homologie d'un \'eclatement}

La d\'emonstration du lemme de d\'eplacement suivant ne pose pas de difficult\'e particuli\`ere.
\begin{prop}\label{movinglemma}
Soit $X$ un sch\'ema r\'egulier int\`egre dans Ob($\cal C$). Soit $Y \subset X$
un sous-sch\'ema ferm\'e propre. Soit $U: X \setminus Y$. 
Pour tout entier $q \geq 0$, l'application
$\oplus_{x \in X_{q}\cap U} \Z  \to CH_{q}(X)$ est surjective.
\end{prop}

La proposition suivante n'est pas non plus trop surprenante.
Elle d\'ecrit le comportement du groupe $H_{2}(\bullet, \Lambda)$
dans un \'eclatement.

\begin{prop}\label{blowupformula}
Soit $X\in Ob({\cal C})$ r\'egulier de dimension $d+1$. Soient $\pi : \tilde{X} \to X$
l'\'eclat\'e de $X$ en un point ferm\'e $x \in X_{s}$, et $E \hookrightarrow  \tilde{X}$
le diviseur exceptionnel.
Supposons $cd_{l}(F) \leq 1$.  Soit $W \subset E \simeq \P^d_{\kappa(x)}$ un sous-sch\'ema ferm\'e int\`egre
de dimension 1.
Si $\Lambda=\Q_{l}/\Z_{l}$, ou si $\Lambda=\Z/l^n$
et le degr\'e 
sur $\kappa(x)$ de $W \hookrightarrow \P^d_{\kappa(x)}$
est premier \`a $l$, alors 
$$ {\rm Ker} \, [ \pi_{*} : H_{2}(\tilde{X},\Lambda) \to H_{2}(X,\Lambda)]= \Lambda(W).$$
 \end{prop}

Pour $\Lambda=\Z/l^n$, le groupe $ \Lambda(W)$ est par d\'efinition 
le sous-groupe de $H_{2}(\tilde{X},\Z/l^n)$ engendr\'e par $\rho_{X}([W])$.
Pour $\Lambda=\Q_{l}/\Z_{l}$, c'est par d\'efinition la limite inductive
de ces groupes pour $n$ tendant vers l'infini.

Pour \'etablir la proposition, on utilise les isomorphismes $$H_{2}(E,\Z/l^n) \simeq H^{2d-2}(E,\Z/l^n(d-1))
\simeq \Z/l^n,$$ o\`u le second r\'esulte de l'hypoth\`ese 
 $cd_{l}(F) \leq 1$ et du calcul bien connu de la cohomologie d'un espace projectif.

\subsection{D\'emonstration du th\'eor\`eme principal}

\begin{theo}\label{surjectif}
Supposons  $R$ hens\'elien et $F$ fini ou s\'eparablement clos. 
Soit $X$ un $R$-sch\'ema dans $\cal{QSP}$.
Pour
 $\Lambda=\Z/l^n$ ou $\Lambda=\Q_{l}/\Z_{l}$,
l'application cycle $$\rho_{X} : CH_{1}(X) \otimes \Lambda \to H_{2}(X,\Lambda)$$
est surjective.
\end{theo}
\noindent{\sc Preuve}   ---  On peut supposer $X$ int\`egre et $\Lambda=\Z/l^n$.
Pour $\dim(X)=1$,  resp. $\dim(X)=2$, la 
  proposition \ref{immediat}, resp.
  \ref{dimension2}, montre que l'application cycle est un isomorphisme.
  Supposons $\dim(X) \geq 3$.  Pour \'etablir l'\'enonc\'e,
 on proc\`ede  par r\'ecurrence sur la dimension. Par  le th\'eor\`eme de
 Bertini relatif, version 1,
  il existe  un diviseur ample $Y \subset X$ qui d\'efinit une $\cal{QSP}$-paire ample.
  Soit $U:=X \setminus Y$. Par la fonctorialit\'e de l'homologie et de l'application
  cycle, on a un diagramme commutatif  
  \[
\xymatrix{
CH_{1}(Y)/l^n \ar[r]  \ar[d]_{\rho_{Y}} &  CH_{1}(X)/l^n  \ar[d]_{\rho_{X}}  & \\
H_{2}(Y,\Z/l^n ) \ar[r]  & H_{2}(X,\Z/l^n )  \ar[r]  & H_{2}(U,\Z/l^n )
}
\]
o\`u la suite inf\'erieure est exacte.  Par hypoth\`ese de r\'ecurrence,
$\rho_{Y}$ est surjective.  Par le th\'eor\`eme  \ref{Lefschetzaffine}  (Lefschetz affine affin\'e)
pour la $\cal{QSP}$-paire ample sur l'anneau de valuation discr\`ete hens\'elien $R$,
comme la dimension de $X$ est au moins 3, on a $ H_{2}(U,\Z/l^n )=0$. Ainsi $\rho_{X}$ est surjectif.
\cqfd

\begin{theo}\label{bijectifinfinidim3}
Supposons  $R$ hens\'elien excellent et $F$ fini ou s\'eparablement clos. 
Soit $X$ un $R$-sch\'ema dans $\cal{QSP}$ de dimension $3$.
L'application cycle $$\rho_{X} : CH_{1}(X) \otimes \Q_{l}/\Z_{l} \to H_{2}(X,\Q_{l}/\Z_{l})$$
est bijective.
\end{theo}
\noindent{\sc Preuve}   ---   
  On a d\'ej\`a \'etabli la surjectivit\'e. Soit $\alpha$ dans le noyau de $\rho_{X}$.
Par le lemme de d\'eplacement \ref{movinglemma}, on peut repr\'esenter $\alpha$ sous la forme
$\sum_{1 \leq j \leq m}[C_{j}] \otimes \lambda_{j}$ avec \mbox{$\lambda_{j} \in \Q_{l}/\Z_{l}$} et
 des courbes ferm\'ees int\`egres   $C_{j} \subset X$   non contenues dans la fibre sp\'eciale $X_{s}$.
Par une version de la r\'esolution plong\'ee des singularit\'es des courbes due \`a U.~Jannsen (appendice de \cite{SS2}, qui utilise l'hypoth\`ese  que $X$ est excellent),
on peut trouver
un compos\'e de $N$~\'eclatements $\pi : \tilde{X} \to X$ en des points ferm\'es de la fibre sp\'eciale
de fa\c con que les transform\'es propres $\tilde{C}_{j}$ des $C_{j}$ soient des courbes r\'eguli\`eres
sans point commun deux \`a deux et satisfassent des conditions de transversalit\'e  convenables
par rapport aux composantes r\'eduites de la fibre sp\'eciale $\tilde{X}_{s}$.

Notons $\tilde{E}_{i} \subset \tilde{X}$, $i=1,\dots, N$ 
les transform\'es propres dans $ \tilde{X}$ des divers diviseurs exceptionnels.

Par le th\'eor\`eme de Bertini relatif, version 2,
il existe
une $\cal{QSP}$-paire ample $(\tilde{X},Y)$ telle que le diviseur ample $\iota : Y \subset \tilde{X}$ contienne
tous les $\tilde{C}_{j}$. On a $\alpha=\pi_{*}\iota_{*}\beta$,
avec $\beta= \sum_{j} [\tilde{C}_{j}] \otimes \lambda_{j} \in CH_{1}(Y) \otimes \Q_{l}/\Z_{l }$.
Comme $Y$ est ample,
 pour chaque $i$,  le sous-sch\'ema $W_{i}:= \tilde{E}_{i} \cap Y$
du sch\'ema  $\tilde{X} $ (de dimension 3)
est  un  ferm\'e non vide, 
int\`egre 
et de dimension 1, et $W_{i}$ n'est inclus dans aucun $\tilde{E}_{i'}$
pour $i' \neq i$. Par une version it\'er\'ee de la proposition  \ref{blowupformula} (homologie d'un \'eclatement),  on a
$$ {\rm Ker} \,  [\pi_{*} : H_{2}(\tilde{X},\Q_{l}/\Z_{l }) \to H_{2}(X,\Q_{l}/\Z_{l }) ] = \sum_{i=1}^N \Q_{l}/\Z_{l }(W_{i}).$$
Dans  le diagramme
   \[
\xymatrix{
  &  & CH_{1}(Y) \otimes \Q_{l}/\Z_{l } \ar[r]^{\iota_{*}} \ar[d]_{\rho_{Y}} &
                  CH_{1}(\tilde{X}) \otimes \Q_{l}/\Z_{l }   \ar[d]^{\rho_{\tilde{X}}} &   \\
                  & H_{3}(\tilde{X} \setminus Y,\Q_{l}/\Z_{l }) \ar[r] & H_{2}(Y,\Q_{l}/\Z_{l }) \ar[r]^{\iota_{*}}  & 
                  H_{2}(\tilde{X},\Q_{l}/\Z_{l })    &  }
\]
la ligne m\'ediane est la suite exacte de localisation, et le carr\'e est commutatif.
Le th\'eor\`eme \ref{Lefschetzaffine}  (th\'eor\`eme de Lefschetz affine affin\'e) assure la
nullit\'e du terme $H_{3}$. Comme $Y$ est dans $Ob({\cal C})$ et de dimension 2,
on sait que la fl\`eche $\rho_{Y}$ est injective.

Combinant ces r\'esultats avec le diagramme commutatif
   \[
\xymatrix{  CH_{1}(\tilde{X}) \otimes \Q_{l}/\Z_{l }   \ar[r]^{\rho_{\tilde{X}}}  \ar[d]^{\pi_{*}} &  H_{2}(\tilde{X},\Q_{l}/\Z_{l })  \ar[d]^{\pi_{*}}\\
   CH_{1}(X) \otimes \Q_{l}/\Z_{l }   \ar[r]^{\rho_{X}} &  H_{2}(X,\Q_{l}/\Z_{l }) 
}
\]
on voit que $ \beta $ est dans le groupe engendr\'e par les $W_{i}$ dans $CH_{1}(Y) \otimes \Q_{l}/\Z_{l } $.
Comme l'image par $\pi : \tilde{X} \to X$ de chaque $W_{i}$ est un point ferm\'e, on conclut $$\alpha=\pi_{*}\iota_{*}(\beta)=0
\in  CH_{1}(X) \otimes \Q_{l}/\Z_{l }.$$
\cqfd

\begin{theo}\label{bijectifdim3}
 Supposons  $R$ hens\'elien excellent et $F$ fini ou s\'eparablement clos. 
Soit $X$ un $R$-sch\'ema dans $\cal{QSP}$ de dimension 3. 
   
{\rm (a)}  L'application cycle $$\rho_{X} : CH_{1}(X) \otimes \Z/l^n \to H_{2}(X,\Z/l^n)$$
est bijective.

{\rm (b)} Pour tout $n>0$ on a  $KH_{3}(X,\Z/l^n) = KH_{3}(X,\Q_{l}/\Z_{l })[l^n]=0.$
\end{theo}
 \noindent{\sc Preuve}   ---   La surjectivit\'e  dans (a) a \'et\'e \'etablie au th\'eor\`eme \ref{surjectif}. 
  Il reste \`a \'etablir l'injectivit\'e.  D'apr\`es la proposition \ref{immediat} (dont la preuve utilise le
 th\'eor\`eme de Merkur'ev et Suslin), pour \'etablir (a) et (b)  
 il suffit de montrer $KH_{3}(X,\Q_{l}/\Z_{l })=0$. On fixe un diviseur $Y \subset X$
 tel que $(X,Y;U)$ soit une $\cal{QSP}$-paire ample. Comme $X$ est de dimension 3,
 la restriction $KH_{3}(X,\Q_{l}/\Z_{l }) \to KH_{3}(U,\Q_{l}/\Z_{l })$ est injective (c'est purement formel, ceci vaut pour tout
 ouvert $U$ dense dans $X$). Il suffit donc de montrer $ KH_{3}(U,\Q_{l}/\Z_{l })=0$.
On a le diagramme commutatif de suites de localisation
    \[
\xymatrix{ 
 CH_{1}(Y) \otimes \Q_{l}/\Z_{l }   \ar[r]  \ar[d]^{\rho_{Y}}   
 &  CH_{1}(X) \otimes \Q_{l}/\Z_{l }   \ar[r] \ar[d]^{\rho_{X}}  
 &  CH_{1}(U) \otimes \Q_{l}/\Z_{l }   \ar[r]  \ar[d]^{\rho_{U}}
   & 0 \\
 H_{2}(Y,\Q_{l}/\Z_{l }) \ar[r]  
 & H_{2}(X,\Q_{l}/\Z_{l })  \ar[r]  
 & H_{2}(U,\Q_{l}/\Z_{l }). 
 &   
 }
\]
Nous avons d\'ej\`a \'etabli que $\rho_{Y}$ et $\rho_{X}$ sont bijectives.
Ainsi $\rho_{U}$ est injective. 
Mais $ H_{2}(U,\Q_{l}/\Z_{l })=0$ (Thm. \ref{Lefschetzaffine}, Lefschetz affine affin\'e).
Donc $ CH_{1}(U) \otimes \Q_{l}/\Z_{l } =0$. Par ailleurs $H_{3}(U,\Q_{l}/\Z_{l }) =0$
(Thm. \ref{Lefschetzaffine}, Lefschetz affine affin\'e). En appliquant la proposition \ref{immediat}(3)   \`a $U$,
 on conclut
 $ KH_{3}(U,\Q_{l}/\Z_{l })=0$.\cqfd
 
 Le th\'eor\`eme principal en dimension quelconque s'\'enonce :
 \begin{theo}\label{thmprincipal}
 Supposons  $R$ hens\'elien excellent et $F$ fini ou s\'eparablement clos. 
Soit $X$ un $R$-sch\'ema dans $\cal{QSP}$.   
 L'application cycle $$\rho_{X} : CH_{1}(X) \otimes \Z/l^n \to H_{2}(X,\Z/l^n)$$
est un isomorphisme de groupes finis.
\end{theo}
  \noindent{\sc Preuve}   ---  La surjectivit\'e  a \'et\'e \'etablie (Thm. \ref{surjectif}).  Pour \'etablir l'injectivit\'e, on fait   une r\'ecurrence sur 
 la dimension,  dans le m\^eme esprit que celle du th\'eor\`eme \ref{bijectifinfinidim3},   avec $\Lambda=\Z/l^n$,
 en utilisant le th\'eor\`eme de Bertini relatif, version 2. Je renvoie ici  le lecteur \`a \cite{SS2}.
 Des isomorphismes
 $H_{2}(X,\Z/l^n)  \simeq H^{2d}(X,\Z/l^n(d))$ (Proposition \ref{dualite})  
et $H^{2d}(X,\Z/l^n(d)) \simeq H^{2d}(X_{s},\Z/l^n(d))$ (changement de base propre, $R$ est hens\'elien)
on d\'eduit que pour $F$ s\'eparablement clos ou fini les groupes $H_{2}(X,\Z/l^n)$ sont finis.
 \cqfd

 \begin{rema}
K. Sato \cite{Sato2} a d\'efini des applications cycle $\rho_{X,p^r}$  sur les quotients $CH_{1}(X)/p^r$,
\`a valeurs dans la cohomologie \'etale de certains complexes de $\Z/p^r$-faisceaux \'etales.
Lorsque $X/R$ est  propre, r\'egulier et semistable, Saito et Sato montrent dans \cite{SS3}
que ces applications sont surjectives.
\end{rema}
 
\subsection{Applications aux groupes de Chow}

 Le th\'eor\`eme suivant est une sorte de th\'eor\`eme de Lefschetz pour le groupe de Chow des 1-cycles
 modulo un entier.
 \begin{theo}\label{Lefschetz1cycles}
 Soit $R$ un anneau de valuation discr\`ete hens\'elien excellent de corps r\'esiduel~$F$.
 Soit $X$ un $R$-sch\'ema dans $\cal{QSP}$.   Soit $n \geq 1$ un entier.
 Soit $i :Y \hookrightarrow X$ un diviseur tel que
  $(X,Y)$ soit une $\cal{QSP}$-paire ample. Soit $U = X \setminus Y$.

{\rm (a)} Supposons $F$ s\'eparablement clos et $\dim(X)=d+1 \geq 2$.  Alors   $CH_{1}(U)/l^n=0$ et
l'application   $i_{*} = CH_{1}(Y)/l^n \to CH_{1}(X)/l^n$ est 
  surjective. Si $d \geq 2$,
elle est bijective.

{\rm (b)} Supposons $F$ fini et $\dim(X)=d+1 \geq 3$.
Alors $CH_{1}(U)/l^n=0$ et
l'application $i_{*} = CH_{1}(Y)/l^n \to CH_{1}(X)/l^n$ est 
  surjective. Si $d \geq 3$,
elle est bijective.

   \end{theo}
   \noindent{\sc Preuve}   --- Consid\'erons le diagramme commutatif \`a lignes exactes
  \[
\xymatrix{ 
&  CH_{1}(Y)/l^n   \ar[r]^{i_{*}}  \ar[d]^{\rho_{Y}}   
 &  CH_{1}(X)/l^n   \ar[r] \ar[d]^{\rho_{X}}  
 &  CH_{1}(U) /l^n   \ar[r]  \ar[d]^{\rho_{U}}
   & 0 \\
H_{3}(U,\Z/l^n) \ar[r] & H_{2}(Y,\Z/l^n) \ar[r]  
 & H_{2}(X,\Z/l^n)  \ar[r]  
 & H_{2}(U,\Z/l^n). 
 &   
 }
\]
D'apr\`es le th\'eor\`eme  \ref{thmprincipal},   $\rho_{Y}$ et $\rho_{X}$ sont des isomorphismes.
Le th\'eor\`eme r\'esulte alors du th\'eor\`eme de Lefschetz affine affin\'e
(Th\'eor\`eme \ref{Lefschetzaffine}). \cqfd

\begin{theo}
Soit $R$ un anneau de valuation discr\`ete hens\'elien excellent de corps r\'esiduel~$F$ s\'eparablement clos.
Soit $X$ un objet de $\cal{QSP}$, et soient $Y_{1}, \dots, Y_{N}$ les composantes de $X_{s,red}$.
Alors pour tout entier $n>0$ premier \`a $car(F)$, l'intersection avec les composantes
induit un isomorphisme
$$ CH_{1}(X)/n \loi  \bigoplus_{i=1}^N \Z/n\Z.$$
\end{theo}
 
\noindent{\sc Preuve}   ---  Soit $d=\dim(X)-1$.
On a des isomorphismes naturels
\begin{align*}
H_{2}(X,\Z/l^n) &\simeq H^{2d}(X,\Z/l^n(d))\\
&\simeq H^{2d}(X_{s},\Z/l^n(d))  \loi  \bigoplus_{i=1}^N  H^{2d}(Y_{j},\Z/l^n(d)) 
 \loi  \bigoplus_{i=1}^N \Z/l^n,
\end{align*}
 et d'apr\`es le th\'eor\`eme principal \ref{thmprincipal} on a un isomorphisme $CH_{1}(X)/l^n \simeq H_{2}(X,\Z/l^n)$. On v\'erifie que pour
 chaque composante $Y_{i}$ la fl\`eche associ\'ee $CH_{1}(X)/l^n \to \Z/l^n$ associe \`a un 1-cycle la classe modulo $l^n$ de
  son nombre d'intersection
 avec le diviseur vertical $Y_{i} \subset X$.
 \cqfd

\begin{theo}\label{bonnereduc}
Soit $R$ un anneau de valuation discr\`ete hens\'elien excellent, de corps des fractions $k$,
de corps r\'esiduel $F$
fini ou s\'eparablement clos. Soit $p$ l'exposant caract\'eristique de $F$.
Soit $X$ un $R$-sch\'ema projectif et lisse \`a fibres g\'eom\'etriquement connexes.
Soit $i : X_{s} \hookrightarrow X$ la fibre sp\'eciale, 
et $j : V =X_{\eta}\hookrightarrow X$  
la fibre g\'en\'erique.

{\rm (1)} Pour tout entier   $n>0$ premier \`a $p={\rm car}(F)$, 
on a des isomorphismes de groupes ab\'eliens finis
$$ CH_{0}(X_{s})/n \lloi  CH_{1}(X)/n \loi CH_{0}(V)/n,$$
o\`u $ CH_{0}(V)$ d\'esigne le groupe de Chow  de dimension z\'ero 
 usuel de la $k$-vari\'et\'e $V$, la fl\`eche de gauche est $i^*$ et la fl\`eche de droite $j^*$.
 
 {\rm (3)} Supposons $F$ s\'eparablement clos. Le groupe $A_{0}(V)$ 
 est divisible par tout entier $n$ premier \`a $p$.
 
{\rm (3)} Supposons $F$ fini. 

(a) Le groupe $A_{0}(V)$ est extension du groupe  fini $ A_{0}(X_{s})$ 
 par un groupe $p'$-divisible. 
 
 (b) Pour tout entier $n$ premier \`a $p$, l'application cycle
 $ CH_{0}(V)/n \to  H^{2d}(V,\mu_{n}^{\otimes d})$ est injective.

(c) Soit $\br'\ (V)$ le sous-groupe de $\br \ V$ de torsion premi\`ere \`a $p$.
Le  noyau  \`a gauche de l'accouplement
$$ A_{0}(V) \times \br '\ V \to \br' k  \subset \Q/\Z$$
est divisible par tout entier $n$ premier \`a $p$.
\end{theo}

\noindent{\sc Preuve}   ---  Soit $d$ la dimension de $V$.
Soit $n$ un entier premier \`a $p$.
Les applications cycle sur $X$ et sur $X_{s}$ 
s'inscrivent dans un diagramme commutatif
\[
\xymatrix{ 
  CH_{1}(X)/n   \ar[r]^{cl_{X}}  \ar[d]^{i^*}
 &  H^{2d}(X,\Z/n(d))   \ar[d]^{i^*} \\
 CH_{0}(X_{s})/n   \ar[r]^{cl_{X_{s}} }  
  &  H^{2d}(X_{s},\Z/n(d)). 
  }
\]
La fl\`eche verticale de droite est un isomorphisme (changement de base propre).
D'apr\`es le th\'eor\`eme principal \ref{thmprincipal},
la fl\`eche $cl_{X}$ est un isomorphisme.
Comme $X$ est lisse sur $R$ et $R$ hens\'elien, la fl\`eche
de restriction $i^* : CH_{1}(X) \to CH_{0}(X_{s})$ est surjective.
Ceci suffit \`a assurer que toutes les fl\`eches dans le diagramme
ci-dessus sont des isomorphismes. Il en est donc ainsi de
$i^* :   CH_{1}(X)/n \to CH_{0}(X_{s})/n$.

 [Lorsque $F$
est un corps fini, on reconna\^{\i}t dans l'isomorphisme
obtenu $ CH_{0}(X_{s})/n  \simeq  H^{2d}(X_{s},\Z/n(d)) $
un cas particulier du th\'eor\`eme  du  corps de classes non ramifi\'e
pour les vari\'et\'es projectives et lisses
sur un corps fini, voir \cite{Sz}.]

Une suite de localisation \'el\'ementaire fournit la suite exacte
 $$CH_{1}(X_{s})/n  \buildrel {i_{*}}  \over \longrightarrow
     CH_{1}(X)/n   \buildrel {j^*}  \over \longrightarrow  CH_{0}(X_{\eta})/n \to 0.$$
Comme $X_{s}$ est un diviseur principal sur $X$, l'application compos\'ee
$$CH_{1}(X_{s})/n  \buildrel {i_{*}}  \over \longrightarrow
     CH_{1}(X)/n   \buildrel {i^{*}}  \over \longrightarrow CH_{0}(X_{s})/n $$
est nulle. Comme on a  \'etabli que  $i^*$ est un isomorphisme, ceci implique que $i_{*}$
est nul. Ainsi $j^*$ est un isomorphisme. Ceci \'etablit le point (1).

Le point (2) en r\'esulte, puisque le groupe $A_{0}(X_{s})$ est divisible.

Supposons $F$ fini. Le th\'eor\`eme de Kato et Saito \cite{KS} assure que le groupe $A_{0}(X_{s})$
est fini. L'application de r\'eduction  $A_{0}(X_{\eta}) \to A_{0}(X_{s})$ est surjective
(lemme de Hensel). En utilisant (1) (voir aussi le 
 th\'eor\`eme \ref{elementaire}) on obtient (3)(a). Comme $X/R$ est lisse,
 une version connue de la conjecture de Gersten assure que l'application
 de restriction $H^{2d}(X,\mu_{n}^{\otimes d}) \to H^{2d}(V,\mu_{n}^{\otimes d})$
 est injective. Du diagramme commutatif
   \[
\xymatrix{ 
  CH_{1}(X)/n   \ar[r]^{j_{*}}  \ar[d]^{cl_{X}}  
 &  CH_{0}(V)/n \ar[d]^{cl_{V}}  
  \\
 H^{2d}(X,\Z/n) \ar[r]  
 & H^{2d}(V,\Z/n),   
 }
\]
 o\`u la fl\`eche sup\'erieure est surjective, on d\'eduit (3)(b).
 L'\'enonc\'e (3)(c) r\'esulte alors du fait connu que l'accouplement naturel
  $$H^{2d}(V,\mu_{n}^{\otimes d}) \times H^2(V,\mu_{n}^{\otimes d}) \to H^{2d+2}(V,\mu_{n}^{\otimes d+1} ) \simeq \Z/n$$
  est une dualit\'e parfaite de groupes finis, et que cet accouplement induit un accouplement
  $$CH_{0}(V)/n \times \br(V) [n] \to \Z/n.$$
 \cqfd

\begin{lemm}\label{architrivial}
Soit  $A$ un groupe ab\'elien.

{\rm (i)} Les propri\'et\'es suivantes sont \'equivalentes

\quad {\rm (a)}  Le groupe $A$ est la somme directe d'un groupe fini d'ordre premier \`a $p$ et d'un groupe
$p'$-divisible.

\quad {\rm (b)}  Pour presque tout premier $l \neq p$, le quotient $A/l$ est nul, et pour tout premier
$l \neq p$, il existe un entier $n_{l}>0$ tel que  la projection $A/l^{n+1} \to A/l^n$ soit
un isomorphisme de groupes finis pour tout $n \geq n_{l}$.

{\rm (ii)} Si un groupe $A$ poss\`ede ces propri\'et\'es, il en est de m\^eme de tout quotient de $A$.\break\hfill \cqfd
\end{lemm}

\begin{theo}\label{applicsemistable}
Soit $R$ un anneau de valuation discr\`ete hens\'elien excellent, de corps des fractions $k$,
de corps r\'esiduel $F$
fini ou s\'eparablement clos. Soit $p$ l'exposant caract\'eristique de $F$.  Soit $V$ une vari\'et\'e projective, lisse, g\'eom\'etriquement connexe sur un  $k$.
 Soit $X$  un $R$-sch\'ema dans $\cal{QSP}$
et soit $V$ la $k$-vari\'et\'e $X_{\eta}$, suppos\'ee g\'eom\'etriquement int\`egre. 
Le groupe $A_{0}(V)$ est isomorphe \`a la somme directe d'un groupe fini
d'ordre premier \`a $p$ et d'un groupe $p'$-divisible.
En particulier $A_{0}(V)/l=0$ pour presque tout premier $l$,
et $A_{0}(V)/n$ est fini pour tout $n>0$ premier \`a $p$.
\end{theo}
\noindent{\sc Preuve}   ---
D'apr\`es
 le th\'eor\`eme \ref{Lefschetz1cycles}
il existe une $R$-courbe $Y \hookrightarrow X$ qui est dans  $\cal{QSP}$
et pour laquelle  pour tout $l$ premier, $l \neq p$ et 
pour tout $n>0$  l'application $CH_{1}(Y)/l^n \to CH_{1}(X)/l^n$ est surjective.
Les restrictions \`a la fibre g\'en\'erique $CH_{1}(X) \to CH_{0}(X_{\eta})$
et $CH_{1}(Y) \to CH_{0}(Y_{\eta})$ sont clairement surjectives, et elles sont compatibles.
Ainsi  les applications naturelles
$CH_{0}(Y_{\eta})/l^n \to CH_{0}(X_{\eta})/l^n$ sont surjectives.
Cela implique le m\^eme \'enonc\'e pour $A_{0}(Y_{\eta})/l^n \to A_{0}(X_{\eta})/l^n$.
Pour la courbe projective et lisse $Y_{\eta}$, on montre 
que
le groupe $A_{0}(Y_{\eta})$ est la somme directe d'un groupe fini
d'ordre premier \`a $p$ 
et d'un groupe  $p'$-divisible.
Le lemme \ref{architrivial} donne alors l'\'enonc\'e sur 
 $A_{0}(X_{\eta})$.
\cqfd
\begin{rema}
La finitude de $A_{0}(V)/n$ et de $CH_{0}(V)/n$ r\'esulte directement de la finitude
de $CH_{1}(X)/n \simeq H_{2}(X,\Z/l^n)$.
\end{rema}

\begin{theo}\label{applicgeneral}
Soit $R$ un anneau de valuation discr\`ete hens\'elien excellent, de corps des fractions $k$,
de corps r\'esiduel $F$
fini ou s\'eparablement clos. Soit $p$ l'exposant caract\'eristique de $F$.  Soit $V$ une vari\'et\'e projective, lisse, g\'eom\'etriquement connexe sur~$k$.
Le groupe $A_{0}(V)$ est isomorphe \`a la somme directe d'un groupe fini
d'ordre premier \`a $p$ et d'un groupe $p'$-divisible.
En particulier $A_{0}(V)/l=0$ pour presque tout premier $l$,
et $A_{0}(V)/n$ est fini pour tout $n>0$ premier \`a $p$.
\end{theo}
\noindent{\sc Preuve}   --- 
Le th\'eor\`eme d'uniformisation de de Jong \cite{dJ}, dans la version raffin\'ee de Gabber,  
th\'eor\`eme d\'ecrit par Illusie  dans \cite{Ill},
implique que, pour tout premier $l\neq p$,  il existe
un $k$-morphisme  propre $p : V' \to V$, g\'en\'eriquement fini de degr\'e $d$ {\it premier \`a} $l$,
et une extension finie d'anneaux de valuation discr\`ete $R'/R$ 
tels que  la  vari\'et\'e $V'$ soit lisse et g\'eom\'etriquement int\`egre sur le corps des fractions de $R'$
et  admette un mod\`ele $X'/R'$ qui soit $\cal{QSP}$. Les propri\'et\'es usuelles
des groupes de Chow  des vari\'et\'es lisses impliquent que le compos\'e
$$ A_{0}(V)  {\buildrel {p^*} \over \longrightarrow} A_{0}(V') {\buildrel {p_{*}} \over \longrightarrow} A_{0}(V)$$
est la multiplication par $d$.
Le th\'eor\`eme r\'esulte alors du th\'eor\`eme \ref{applicsemistable}.
\cqfd

\begin{rema}
Le th\'eor\`eme originel de de Jong \cite{dJ} combin\'e avec le th\'eor\`eme  \ref{applicsemistable}
suffit \`a \'etablir la trivialit\'e de $A_{0}(V)/l$ pour presque tout premier $l$. C'est cet \'enonc\'e
que l'on trouve dans \cite{SS2}.
\end{rema}

\begin{coro}\label{finitude}
Soient $R$, $p$ et $k$ comme dans le th\'eor\`eme ci-dessus.
Soit $V$ une vari\'et\'e projective, lisse, g\'eom\'etriquement connexe sur  $k$. 
Si pour tout corps alg\'ebriquemen clos $\Omega$ contenant $k$ on a $A_{0}(V\times_{k}\Omega)=0$,
alors le groupe $A_{0}(V)$ est la somme directe d'un groupe fini et d'un groupe
d'exposant une puissance de $p$.
\end{coro}

\noindent{\sc Preuve}   ---   Sur un corps quelconque, on montre  en effet que le groupe $A_{0}(V)$ de toute telle vari\'et\'e $V$
est annul\'e par un entier positif. L'\'enonc\'e  du corollaire r\'esulte alors du th\'eor\`eme. \cqfd

Le r\'esultat s'applique en particulier aux $k$-vari\'et\'es
g\'eom\'etriquement rationnellement connexes (au sens de Koll\'ar, Miyaoka, Mori).
Ainsi le groupe $A_{0}(V)$ d'une vari\'et\'e rationnellement connexe sur le corps $k=\C((t))$
est un groupe fini. C'est une question ouverte de savoir si dans ce cas le
groupe $A_{0}(V)$ est nul.

\subsection{Applications aux conjectures  \ref{conjsepclos} et \ref{conjfini} }

\begin{theo}\label{annulationKH2}
Soit $R$ un anneau de valuation discr\`ete hens\'elien excellent de corps r\'esiduel~$F$.
 Soit $X$ un $R$-sch\'ema dans $\cal{QSP}$. Soit $n \geq 1$ un entier.

{\rm (a)} Si $F$ est un  corps s\'eparablement clos,   alors  $K\!H_{2}(X,\Z/l^n)\!=\!0$ 
et \, \mbox{$K\!H_{2}(X,\Q_{l}/\Z_{l })\!=\!0$.}

{\rm (b)} Si $F$ est un corps fini, on a $KH_{2}(X,\Q_{l}/\Z_{l })=0$.
\end{theo}
\noindent{\sc Preuve}   ---  Dans chacun des cas consid\'er\'es,
la proposition \ref{annulationgaloisienne} assure que la suite spectrale
de niveau est concentr\'ee dans le premier quadrant. 
La forme de cette suite spectrale donne alors
des suites exactes
$$ CH_{1}(X) \otimes \Lambda \to H_{2}(X,\Lambda) \to KH_{2}(X,\Lambda) \to 0.$$
Du th\'eor\`eme \ref{surjectif} on d\'eduit $KH_{2}(X,\Lambda)=0$.
\cqfd

\begin{theo}\label{annulationKH3}
Soit $R$ un anneau de valuation discr\`ete hens\'elien excellent de corps r\'esiduel~$F$.
 Soit $X$ un $R$-sch\'ema dans $\cal{QSP}$. Soit $n \geq 1$ un entier.

{\rm (a)} Si $F$ est un  corps s\'eparablement clos, alors  $K\!H_{3}(X,\Z/l^n)\!=\!0$ 
et  \mbox{$K\!H_{3}(X,\Q_{l}/\Z_{l })\!=\!0$.}

{\rm (b)} Si $F$ est un corps fini et  $\dim(X) \leq 4$, alors $KH_{3}(X,\Z/l^n)=0$.

{\rm (c)} Si $F$ est un corps fini, on a $KH_{3}(X,\Q_{l}/\Z_{l })=0$.
\end{theo}
\noindent{\sc Preuve}   --- 
L'\'enonc\'e est  trivial pour $\dim(X) \leq 2$, et il a \'et\'e d\'emontr\'e pour
$\dim(X)=3$ (Thm. \ref{bijectifdim3}).
Supposons
donc $\dim(X) \geq 4$ et \'etablissons le r\'esultat par r\'ecurrence sur la dimension.
Soit $(X,Y;U)$ une $\cal{QSP}$-paire ample. 
Pour  $\Lambda=\Z/l^n$ si $F$ est s\'eparablement clos
 et $\Lambda=\Q_{l}/\Z_{l }$ si $F$ est fini,
pour tout $X \in {\cal C}$,  la proposition \ref{annulationgaloisienne} 
assure que la suite spectrale
de niveau est concentr\'ee dans le premier quadrant
(c'est ici qu'on se limite \`a $\Lambda=\Q_{l}/\Z_{l }$
lorsque $F$ est un corps fini).
Appliquant ceci \`a $U$,
on trouve une suite exacte 
$$ H_{3}(U,\Lambda) \to KH_{3}(U,\Lambda) \to CH_{1}(U)\otimes \Lambda.$$

Lorsque $\dim(U)=4$, la forme de la suite spectrale assure que l'on a encore
cette suite exacte pour $F$ fini et $\Lambda=\Z/l^n$.

D'apr\`es le th\'eor\`eme \ref{Lefschetz1cycles}, on a $CH_{1}(U)\otimes \Lambda=0$.
D'apr\`es le th\'eor\`eme de Lefschetz affine affin\'e \ref{Lefschetzaffine}, 
comme on a $\dim(U) \geq 4$, on a $H_{3}(U,\Lambda)=0$. On conclut donc
$KH_{3}(U,\Lambda)=0$ dans chacun des trois cas (a), (b), (c).

Par ailleurs on a une longue suite exacte
$$ \dots  \to KH_{3}(Y,\Lambda) \to KH_{3}(X,\Lambda) \to KH_{3}(U,\Lambda) \to \dots $$
Par hypoth\`ese de r\'ecurrence, $ KH_{3}(Y,\Lambda)=0$. Ainsi $ KH_{3}(X,\Lambda)=0$.
\cqfd

\begin{rema}
Dans \cite{SS2}, les auteurs demandent si l'\'enonc\'e (b) vaut en toute dimension.
 En \'etablissant le th\'eor\`eme \ref{bijectifdim3}, 
ils montrent que c'est le cas pour $\dim(X)=3$.
 Dans l'argument ci-dessus,
ils commencent la r\'ecurrence
 en dimension 3, et 
n'observent pas le r\'esultat pour $\dim(X)=4$.
 \end{rema}

\section{R\'esultats r\'ecents et questions ouvertes}

Le th\'eor\`eme de Saito et Sato joue un r\^ole important dans la d\'emonstration
du th\'eor\`eme suivant, dont la d\'emonstration, qui utilise la th\'eorie de Hodge et est tr\`es \'elabor\'ee,
ne peut \^etre \'evoqu\'ee ici.
 
\begin{theo}[M. Asakura et S. Saito \cite{AS}]
\label{asakurasaito}
Soient $k$ un corps $p$-adique, $R$ son anneau d'entiers et $F$ son corps r\'esiduel.
Soit $X \subset \P^3_{R}$ une $R$-hypersurface  lisse de degr\'e au moins 5.
Supposons la fibre g\'en\'erique $X_{k}$ tr\`es g\'en\'erale.
Soit $r$ le rang du groupe de Picard
de la fibre sp\'eciale $X_{F}$. Alors le sous-groupe
de torsion $l$-primaire de $A_{0}(X)$ est somme d'un groupe fini
et de $(\Q_{l}/\Z_{l})^{r-1} $.
\end{theo}
Il est facile de donner des exemples de telles surfaces avec $r>1$.

\bigskip

 Le th\'eor\`eme suivant g\'en\'eralise une partie du th\'eor\`eme \ref{RoquetteLichtenbaum},
qui porte sur les courbes, pour lesquelles $\br X=0$ (Th\'eor\`eme \ref{ArtinGrothendieck}).

\begin{theo}[S. Saito et K. Sato \cite{SS3}]
Soient $k$ un corps $p$-adique, $R$ son anneau d'entiers, $X$ un $R$-sch\'ema
propre, connexe et r\'egulier. Le sous-groupe $\br X \subset \br X_{k}$
est dans le noyau \`a droite de l'accouplement
$$ CH_{0}(X_{k}) \times \br X_{k} \to \Q/\Z.$$
 Si le th\'eor\`eme de puret\'e vaut pour le groupe de Brauer de $X$,
alors  $\br X \subset \br X_{k}$ est le noyau \`a droite de cet accouplement.
\end{theo}

La partie premi\`ere \`a $p$ de ce th\'eor\`eme est d\'ej\`a dans \cite{CTS}.
La partie $p$-primaire est  beaucoup plus d\'elicate.

Le th\'eor\`eme de puret\'e vaut pour le groupe de Brauer
 si la dimension de $X$ est au plus 3 (Gabber),
il vaut pour la torsion non $p$-primaire du groupe
de Brauer  (Gabber \cite{Fuj}),   
il vaut aussi pour la torsion $p$-primaire dans un certain nombre de cas (voir \cite{SS3}).

\bigskip

{\it Quelques questions}

\medskip

Soient $k$ un corps $p$-adique et $V$ une $k$-vari\'et\'e projective,  lisse, g\'eom\'etriquement connexe.

(1) Pour $V$ de dimension au moins 3 et $n>0$ entier,
  le groupe $A_{0}(V)[n]$ est-il fini ?
  
(2) Le quotient $A_{0}(V)/p$ est-il fini ?

(3) Les  noyaux des applications $\alb_{X} :  A_{0}(V) \to \Alb_{V}(k)$
et $A_{0}(V) \to \Hom(\br V, \Q/\Z)$  sont-ils chacun extension d'un groupe fini
par un groupe divisible ?

(4) Supposons que $V/k$ est la fibre g\'en\'erique de 
$X/R$ projectif quasisemistable sur l'anneau des entiers $R$ de $k$. Pour $n$ entier premier \`a $p$,
le th\'eor\`eme principal \ref{thmprincipal} donne une formule pour le quotient $CH_{1}(X)/n$.
Peut-on en d\'eduire une formule pour $CH_{0}(V)/n$ ?
Un cas particulier est \'etudi\'e dans \cite{Da}.

(5) Peut-on comprendre de fa\c con \og invariante \fg \ l'exemple de Parimala et Suresh \cite{PS} ?

\end{document}